\input amstex
\documentstyle{amsppt}
\magnification 1200
\vcorrection{-1cm}
\NoBlackBoxes
\input epsf

\rightheadtext{Algorithmic recognition of quasipositive braids}
\topmatter
\title     Algorithmic recognition of quasipositive 4-braids of algebraic length three
\endtitle
\author   S.Yu.~Orevkov
\endauthor
\address  IMT, Universit\'e Toulouse-3, France \endaddress
\address Steklov Math. Institute, Moscow, Russia \endaddress
\abstract
We give an algorithm to decide whether a given braid with four strings
is a product of three factors which are conjugates of standard generators
of the braid group. The algorithm is of polynomial time. It is based on the Garside theory.
We give also a polynomial algorithm to decide if a given braid with any number of strings
is a product of two factors which are conjugates of given powers of the standard generators
(in my previous paper this problem was solved without polynomial estimates).
\endabstract
\endtopmatter

\def\SSS{\operatorname{SSS}}
\def\USS{\operatorname{USS}}
\def\SC{\operatorname{SC}}

\def\SCRev{\operatorname{SC}^{\Lsh}}
\def\veeRev{\vee^{\Lsh}}
\def\wedgeRev{\wedge^{\!\Lsh}}
\def\cRev{\bold c^{\!\Lsh}}

\def\dRev{\bold d^{\Lsh}}
\def\sRev{\frak s^{\Lsh}}
\def\iotaRev{\iota^{\Lsh}}
\def\phiRev{\varphi^{\Lsh}}
\def\varphiRev{\varphi^{\Lsh}}

\def\SCRev{{\SC}^{\Lsh}}
\def\Br {\operatorname{Br}}
\def\Z{\Bbb Z}

\def\refBessis         {1}
\def\refBessisCorran   {2} \def\refBC{\refBessisCorran}
\def\refBGGMi          {3} 
\def\refBKLone         {4}
\def\refBKLtwo         {5}
\def\refBO             {6}
\def\refCW             {7}        
\def\refCharney        {8}
\def\refDehornoy       {9}
\def\refDP             {10}       
\def\refEM             {11}       
\def\refGarside        {12}
\def\refGebGM          {13} \def\refGG{\refGebGM}
\def\refLee            {14}
\def\refOrevkovTop     {15}
\def\refOrevkovUR      {16}
\def\refOrevkovGAFA    {17}
\def\refOrevkovQPthree {18}
\def\refOrW            {19}
\def\refQPtwo          {20}
\def\refPicant         {21}
\def\refRudolph        {22}


\def\sectGarside {2}
\def\sectDef     {\sectGarside.1}
\def\sectPrelim  {\sectGarside.2}

\def\sectSSS     {3}

\def\sectQPfour  {4}

\def\sectProof   {5}

\def\sectLee     {6}


\def\thSSS     {1.1}
\def\corSSS    {1.2}
\def\corUSS    {1.3}
\def\thMain    {1.4}
\def\remMain   {1.5}
\def\remExample{1.6}

\def\lemInfSup   {\sectGarside.1}
\def\lemInf                       {\lemInfSup a}
\def\lemSup                       {\lemInfSup b}
\def\lemXtY      {\sectGarside.2}
\def\lemCharney  {\sectGarside.3}
\def\corCharney  {\sectGarside.4}
\def\lemCharneyR {\sectGarside.5}

\def\lemAx       {\sectSSS.1}
\def\thQPone     {\sectSSS.2}
\def\thBlock     {\sectSSS.3}
\def\lemAAP      {\sectSSS.4}
\def\remAAP      {\sectSSS.5}
\def\lemXt       {\sectSSS.6}
\def\lemUXt      {\sectSSS.7}
\def\propTwoAtoms{\sectSSS.8}
\def\remTwoAtoms {\sectSSS.9}
\def\propLNF     {\sectSSS.10}
\def\lemG        {\sectSSS.11}
\def\lemGp       {\sectSSS.12}
\def\lemGpp      {\sectSSS.13}
\def\propSSSc    {\sectSSS.14}

\def\lemCount {\sectQPfour.1}
\def\lemC     {\sectQPfour.2}
\def\lemD     {\sectQPfour.3}
\def\lemCDi   {\sectQPfour.4}
\def\lemCDii  {\sectQPfour.5}
\def\lemUxVi  {\sectQPfour.6}
\def\lemUxVs  {\sectQPfour.7}
\def\lemForCaseII{\sectQPfour.8}

\def\lemCaseI    {\sectProof.1}
\def\lemCaseII   {\sectProof.2}
\def\lemCaseIII  {\sectProof.3}
\def\lemFinal    {\sectProof.4}

\def\propLee     {\sectLee.1}
\def\lemLeeCC    {\sectLee.2}
\def\lemNonRigid {\sectLee.3}


\def\eqSSSi     {1}
\def\eqSSSii    {2}
\def\eqUXt      {3}
\def\eqUXtI     {4}
\def\eqLemGii   {5}
\def\eqLemGiii  {6}
\def\eqLemGpI   {7}
\def\eqLemGpII  {8}

\def\eqDefCount   {9}
\def\eqLemUxVs    {10}
\def\eqForCaseII  {11}
\def\eqForCaseIIa {12}
\def\eqCaseI     {13}
\def\eqCaseIIinf {14}
\def\eqCaseII    {15}
\def\eqCaseIIa   {16}
\def\eqCaseIIb   {17}
\def\eqCaseIII   {18}
\def\eqCaseIIIa  {19}
\def\eqCaseIIIaa {20}
\def\eqCaseIIIb  {21}
\def\eqCaseIIIc  {22}

\def\figD     {1}
\def\figC     {2}

\document


\head   1. Introduction and statement of main results
\endhead

In this paper we continue the study started in [\refOrevkovQPthree] and [\refQPtwo].
Let $G$ be a Garside group with set of atoms $\Cal A$, for example,
$G=\Br_n$ -- the braid group and $\Cal A=\{\sigma_1,\dots,\sigma_{n-1}\}$ -- the set of its
standard generators (called also Artin generators).
Recall that $\Br_n$ is generated by $\Cal A$
subject to the relations
$$
   \text{$\sigma_i\sigma_j=\sigma_j\sigma_i$ for $|i-j|>1$;}\qquad
   \text{$\sigma_i\sigma_j\sigma_i=\sigma_j\sigma_i\sigma_j$ for
             $|i-j|=1$}.
$$

If an element of $G$ is a product of conjugates of atoms, we say that it
is {\it $\Cal A$-quasipositive}
or just {\it quasipositive} when it is clear which $\Cal A$ is meant.
Note that for Artin-Tits groups (in particular, for braid groups) the notion
of quasipositivity does not depend on the choice between the standard or the dual Garside structure.
We are looking for a solution to the {\it Quasipositivity Problem} -- the algorithmic problem
to decide whether a given element of $G$ is quasipositive or not.
This problem arises in the study of plane complex algebraic or pseudoholomorphic curves,
see, e.~g., [\refRudolph, \refBO, \refOrevkovTop--\refOrevkovGAFA, \refOrW].

Let $e:G\to\Z$ be the homomorphism which takes all atoms to 1.
The value $e(X)$ is called the {\it algebraic length} or {\it exponent sum} of $X$.
The quasipositivity problem for $n$-braids is solved in [\refOrevkovQPthree] for $n=3$ and
in [\refQPtwo] for any $n$ but only for braids of algebraic length two.
Note that the case $n<3$ is trivial and the case $e(X)<2$
is the simplest particular case of the conjugacy problem. 
The case $n=4$, $e(X)=3$ is done in the present paper, see Theorem \thMain.

In fact, a slightly more general problem is solved in [\refQPtwo]. We found an algorithm to decide whether a given
braid $X$ is a product of two conjugates of atom powers. The algorithm in [\refQPtwo]
is rather efficient in practice but no polynomial time bounds are known for it.
Here we give a polynomial time solution to this problem; in the case of braid groups, it is
also polynomial with respect to the number of strings.
Namely, Theorem \thSSS\ states that if $X$ is a product of two conjugates of atom powers,
then each element of the super summit set $\SSS(X)$ for the Birman--Ko--Lee Garside structure
satisfies a certain quickly checkable condition (see Corollary \corSSS\ and Proposition \propLNF), and
it is known [\refBKLtwo] that an element of $\SSS(X)$ can be computed in polynomial time.

Theorem \thSSS\ also plays a central role in our proof of Theorem \thMain\ (the main result of the paper)
which states that if a 4-braid $X$ with $e(X)=3$ is quasipositive, then
$\SSS(X)$ contains an element of the form $xY$ for an atom $x$ and a quasipositive braid $Y$
of algebraic length 2.
So, Theorem \thMain\ solves the quasipositivity problem for $4$-braids $X$ with $e(X)=3$.
This solution is of polynomial time
provided a polynomial upper bound for the size of $\SSS(X)$. Such a bound is given by S.-J.~Lee
[\refLee; Corollary 4.5.4]. 
Note that recently Calvez and Wiest [\refCW] independently obtained
the main result of [\refLee; Chapter 4] (a polynomial time solution to the conjugacy problem in $\Br_4$)
by similar methods.

Let us give precise statements of the main results. For elements $a,b$ of a group $G$ we set
$b^a=a^{-1}ba$, $b^G=\{b^c\mid c\in G\}$, and we write $a\sim b$ if $a\in b^G$.
When speaking of Garside groups,
we use the terminology and notation from [\refQPtwo] which is mostly the same as in [\refGG];
see Section \sectDef\ for a very brief summary.

\proclaim{ Theorem \thSSS }
Let $(G,\Cal P,\delta)$ be a homogeneous symmetric square free Garside structure
of finite type
{\rm(for example, the
Birman-Ko-Lee Garside structure on $\Br_n$)} and let $\Cal A$ be the set of atoms.

Let $Z\in\SSS(Z)\cap\big( (x^k)^G(y^l)^G\big)$ where $k,l\ge1$ and
$x,y\in\Cal A$. 
Then, up to exchange of $x^k$ and $y^l$, one of the following possibilities takes place:
\roster
\item"\rm{(i)}"
$Z=XY$ where $X\sim x^k$, $Y\sim y^l$, and $\ell(Z)=\ell(X)+\ell(Y)$;
\item"\rm{(ii)}"
$Z=x_1^p Y x_1^{k-p}$ where $Y\sim y^l$, $x_1\in x^G\cap\Cal A$,
    $0\le p\le k$, and $\ell(Z) = k+\ell(Y)$;
\item"\rm{(iii)}"
$Z=x_1^p y_1^l x_1^{k-p}$ where $x_1\in x^G\cap\Cal A$, $y_1\in y^G\cap\Cal A$,
    and $0\le p\le k$.
\endroster
\endproclaim


Using the blocking property [\refQPtwo; Corollary 7.2] (see Theorem \thBlock),
Theorem \thSSS\ implies

\proclaim{ Corollary \corSSS } Let the hypothesis of Theorem \thSSS\ hold and $\inf Z<0$.


\smallskip
If Case (i) occurs, i.~e., if $Z=(x_1^k)^P(y_1^l)^Q$ with
$x_1\in x^G\cap\Cal A$,  $y_1\in y^G\cap\Cal A$, and
$\ell(P)+\ell(Q)\ge 1$
{\rm(}we may assume also that $\inf P=\inf Q=0${\rm)}
then the left normal form of $Z$ is
$$
   \delta^{-p-q}\cdot A_1\cdot\dots\cdot A_p
    \cdot C_1\cdot\dots\cdot C_{k+p+q}\cdot y_1^l
    \cdot B_1\cdot\dots\cdot B_q                                         \eqno(\eqSSSi)
$$
where
$A_1\cdot\dots\cdot A_p$,
$\; C_1\cdot\dots\cdot C_{k+p+q}$, and
$B_1\cdot\dots\cdot B_q$ are the left normal forms of
$\delta^p \tau^{-q}(P^{-1})$,
$\delta^q x_1^k P Q^{-1}$, and
$Q$ respectively.


\smallskip
If Case (ii) occurs, i.~e., if $Z=x_1^p (y_1^l)^Q x_1^{k-p}$ with
$x_1\in x^G\cap\Cal A$,  $y_1\in y^G\cap\Cal A$, and $\ell(Q)=n\ge 1$
{\rm(}we may assume also that $\inf Q=0${\rm)}
then the left normal of $Z$ is
$$
   \delta^{-n}\cdot C_1\cdot\dots\cdot C_{p+n}
   \cdot y_1^l\cdot B_1\cdot\dots\cdot B_n\cdot x_1^{k-p}                    \eqno(\eqSSSii)
$$
where
$C_1\cdot\dots\cdot C_{n+p}$ and
$B_1\cdot\dots\cdot B_n$ are the left normal forms of
$\delta^n x_1^p Q^{-1}$ and
$Q$ respectively. \qed
\endproclaim

All possibilities for the left normal forms of $Z$ in Case (iii) of Theorem \thSSS\ are listed in
Proposition \propLNF.

Note that due to Corollary \corSSS, it is very fast to check whether $Z$ satisfies Conditions (i) or (ii):
it is enough to recognize the pattern $y_1^l$ in the left normal form of $Z$ and to check (using Theorem \thQPone)
whether we obtain a conjugate of $x^k$ after its removal; then, of course, the same should be done with
$x^k$ and $y^l$ swapped. If $\inf Z\ge 0$, then Condition (iii) can be checked for all pairs of
atoms $(x_1,y_1)$ from $(x^G)\times(y^G)$ (Proposition \propLNF\ can be used to reduce the number of tests).


\proclaim{ Corollary \corUSS } Let the hypothesis of Theorem \thSSS\ holds and $\inf Z<0$.
Then any cycling orbit of $\USS(Z)$ and any decycling orbit of $\USS(Z^{-1})^{-1}$
contains an element
whose left normal form  is as in [\refQPtwo; Theorem 1b], i.~e., of the form (\eqSSSii) with $p=0$.
\endproclaim

This fact was conjectured in [\refQPtwo; Remark (4) on p.~1083].
In particular, it gives a proof of  [\refQPtwo; Theorem 1b] independent of
the transport properties of cyclic sliding.
Theorem \thSSS\ and Corollary \corUSS\ are proven in Section \sectSSS.
An important ingredient of the proof is
the blocking property of square free homogeneous symmetric Garside structures
[\refQPtwo; Section 7] (see Theorem \thBlock).



\proclaim{ Theorem \thMain }
Let $(G,\Cal P,\delta)$ be a square free homogeneous symmetric Garside structure
of finite type such that $\|\delta\|=3$
{\rm(for example, the Birman-Ko-Lee Garside structure on $\Br_4$\rm)} and let
$\Cal A$ be the set of atoms.

Let $X\in a_1^G a_2^G a_3^G$ with $a_1,a_2,a_3\in\Cal A$.
Then there exists a permutation $(x,y,z)$ of $(a_1,a_2,a_3)$ such that
$\SSS(X)$ contains an element of the form
$x_1 Y$ with $x_1\in x^G\cap\Cal A$, $Y\in y^G z^G$ such that either $\inf Y=\inf x_1 Y$
or $Y\in\Cal P$.
\endproclaim

So, this theorem reduces the quasipositivity problem for the case $e(X)=3$ to the quasipositivity
problem for the case $e(X)=2$.
Theorem \thMain\ is an immediate consequence of Lemmas \lemCaseI\ -- \lemFinal.

\medskip\noindent
{\bf Remark \remMain.} It seems plausible that Theorem \thMain\ holds with minor changes
for products of three conjugates of given powers of atoms.

\medskip\noindent
{\bf Remark \remExample.} The following example shows that $\SSS(X)$ cannot be replaced
by $\USS(X)$ in Theorem \thMain. We consider
the 4-braid $X=\sigma_2^{\sigma_1\sigma_3^3} \sigma_2^{\sigma_1^2\sigma_2^{-1}} \sigma_3^{\sigma_2}$.
Then, for the Birman--Ko--Lee Garside structure on $\Br_4$, we have:
$\ell_s(X)=12$, $\inf_s X=-5$, $\sup_s X=7$, all elements of $\USS(X)$ are rigid, and
$|\USS(X)|=48$. A computation shows that $x^{-1}Z$ is not quasipositive for any $x\in\Cal A$,
$Z\in\USS(X)$.


\medskip
In Section \sectLee\ we give a summary of those results from Lee's thesis [\refLee] about
the structure of $\SSS(X)$ which extend to any homogeneous Garside group with $\|\Delta\|=3$.
This section is independent of the rest of the paper.



\head\sectGarside.   Garside groups
\endhead

\subhead\sectDef. Notation and some definitions
\endsubhead
Given two elements $a,b$ of a group $G$, we set $b^a=a^{-1}ba$ and $b^G=\{b^c\mid c\in G\}$.

Garside groups were introduced in [\refDP, \refDehornoy] as a class of groups to which
Garside's methods [\refGarside] extend.
We use the definitions and notation for Garside structures introduced in [\refGG]
and reproduced almost without changes in [\refQPtwo].
So, a Garside structure on a group $G$ is $(G,\Cal P,\Delta)$ where $\Delta$ is the Garside element
and $\Cal P=\{X\mid X\succcurlyeq 1\}$; we set $\tau(X)=X^\Delta$;
we denote the infimum, supremum, canonical length, and (when $X\in\Cal P$) letter length
of $X\in G$ by $\inf X$, $\sup X$, $\ell(X)$, and $\|X\|$ respectively; we denote the minimal values
of $\inf Y$, $\sup Y$, and $\ell(Y)$ over all $Y\in X^G$ by $\inf_s X$, $\sup_s X$, and $\ell_s(X)$
(see details in [\refGG, \refQPtwo]).

The only difference between the notation in [\refGG] and in [\refQPtwo] is that
we denote the set of simple elements by
$[1,\Delta]$ instead of the commonly used notation $[0,1]$.
We set also
$]1,\Delta]=[1,\Delta]\setminus\{1\}$,
${[1,\Delta[}=[1,\Delta]\setminus\{\Delta\}$,
${]1,\Delta[}={[1,\Delta[}\setminus\{1\}$.

The only new terminology introduced in [\refQPtwo] is the following. We say that
a Garside structure is {\bf homogeneous} if $\|XY\|=\|X\|+\|Y\|$ for any $X,Y\in\Cal P$.
In this case we define a group homomorphism $e:G\to\Bbb Z$ by setting
$e(X)=\|X\|$ for $X\in\Cal P$.
A Garside structure is called {\bf symmetric} if
$A\preccurlyeq B\Leftrightarrow B\succcurlyeq A$ for any simple elements $A,B$ and it is called
{\bf square free} if $x^2\not\preccurlyeq\Delta$ for any atom $x$.
The main example of symmetric homogeneous square free Garside structures are the dual Garside structures
on Artin-Tits groups of spherical type introduced by Bessis [\refBessis], in particular, the
Birman-Ko-Lee Garside structure  [\refBKLone] on $\Br_n$.
Another example is the Garside structure on the braid extension of the complex reflection group
$G(e,e,r)$ introduced in [\refBC].

In this paper we denote the Garside element by $\Delta$ when we speak
of an arbitrary Garside structure, but
we denote it by $\delta$ (as in [\refBKLone]) if
the Garside structure under consideration is supposed to be homogeneous and symmetric.

We denote the left (resp. right) gcd and lcm of $X$ and $Y$ by $X\wedge Y$ and $X\vee Y$
(resp. by $X\wedgeRev Y$ and $X\veeRev Y$).
We denote the usual (i.~e., left) cycling, decycling, and cyclic sliding operators by
$\bold c$, $\bold d$, and $\frak s$ respectively. We denote the initial factor, final factor, and
preferred prefix of $X$ by $\iota(X)$, $\varphi(X)$, and $\frak p(X)$.
So, $\bold c(X)=X^{\iota(X)}$, $\bold d(X)=X^{\varphi(X)^{-1}}$, $\frak s(X)=X^{\frak p(X)}$.
We denote the right counterparts of $\bold c$, $\bold d$, $\iota$, $\varphi$ by
$\cRev$, $\dRev$, $\iotaRev$, $\varphiRev$, i.~e.,
if $A_1\cdot\dots\cdot A_r\cdot \Delta^p$, $r\ge1$,
is the right normal form of $X$, then
$$
\iotaRev(X)=\tau^p(A_r),\quad
\phiRev(X)=A_1,\quad
\cRev(X)=X^{\iotaRev(X)^{-1}},\quad
\dRev(X)=X^{\phiRev(X)}.
$$


\subhead\sectPrelim. Some facts about general Garside groups
\endsubhead
Let $(G,\Cal P,\Delta)$ be any Garside structure of finite type.

\proclaim{ Lemma \lemInfSup } Let $X,Y\in G$. Then:

(a). $\inf XY > \inf X + \inf Y$ if and only if $\Delta\preccurlyeq\iota^\Lsh(X)\iota(Y)$.

(b). $\sup XY < \sup X + \sup Y$ if and only if $\varphi(X)\varphi^\Lsh(Y)\preccurlyeq\Delta$.
\endproclaim

\demo{ Proof }
(a). See [\refQPtwo; Lemma 2.4].

\smallskip
(b). Follows from (a) applied to $Y^{-1}$ and $X^{-1}$.
Indeed, suppose that $\sup XY < \sup X + \sup Y$. Then
$\inf Y^{-1} X^{-1} = \inf(XY)^{-1}=-\sup XY > -\sup X -\sup Y =\inf X^{-1}+\inf Y^{-1}$.
Hence $\Delta\preccurlyeq\iotaRev(Y^{-1})\iota(X^{-1})$ by (a). Note that
$\varphi(X)\iota(X^{-1})=\iotaRev(Y^{-1})\phiRev(Y)=\Delta$, thus
$\Delta\preccurlyeq\iotaRev(Y^{-1})\iota(X^{-1})=
(\Delta\phiRev(Y)^{-1})(\varphi(X)^{-1}\Delta)$ whence
$1\preccurlyeq
\phiRev(Y)^{-1}\varphi(X)^{-1}\Delta$ and, finally,
$\varphi(X)\phiRev(Y)\preccurlyeq\Delta$.
\qed\enddemo

\proclaim{ Lemma \lemXtY } Let $\sup XsY \le\sup X+\sup Y$ where
$X,Y\in G$, $s\in[1,\Delta]$.
Then there exist $u,v\in[1,\Delta]$
such that $s=uv$, $\sup Xu=\sup X$, and $\sup vY=\sup Y$.
\endproclaim

\demo{ Proof }
If $\sup Xs \le\sup X$, then we just set $u=s$, $v=1$ and we are done. So, assume that
$\sup Xs = \sup X+1$. Then, by Lemma \lemSup, we have $\varphi(Xs)\phiRev(Y)\preccurlyeq\Delta$.
Let $v=\varphi(Xs)$. Then $s\succcurlyeq v$ by Lemma \lemCharneyR, i.~e., $s=uv$ for some
$u\in[1,\Delta]$. Since $v=\varphi(Xuv)$, we have $\sup Xuv =\sup Xu +\sup v $, hence
$\sup Xu =\sup Xs -\sup v =\sup Xs - 1=\sup X$. Since $v\phiRev(Y)\preccurlyeq\Delta$,
we have $\sup vY =\sup Y$.
\qed\enddemo

\proclaim{ Lemma \lemCharney }  {\rm[\refCharney; Prop.~3.1].}
Suppose that $X=A_1\cdot A_2\cdot\dots\cdot A_r$
is in left normal form
{\rm(}$A_i\in{]1,\Delta[}$, $i=1,\dots,r${\rm)},
and let $A_0$ be a simple element.
Then the decomposition $A_0X=A_0'\cdot A_1'\cdot\dots\cdot A_r'$ is left weighted
where the $A'_i$'s are defined recursively together with simple elements $t_0,\dots,t_r$ by
the conditions that
 $t_0=A_0$,  $A'_{i-1}\cdot t_i$ is the left normal form of $t_{i-1}A_i$ for $i=1,\dots,r$, and $A'_r=t_r$.
We have  $A'_i\ne\Delta$ for $i>0$ and $A'_i\ne1$ for $i<r$ {\rm(}but
it is possible that $A'_0=\Delta$ or $A'_r=1${\rm)}.
\qed\endproclaim

\proclaim{ Corollary \corCharney } Under the hypothesis of Lemma \lemCharney, suppose that
$\sup A_0X=\sup A_0 + \sup X$ and $\|A_i\|=1$ for some $i\in\{1,\dots,r\}$.
Then $\varphi(A_0X)=\varphi(X)$. \qed
\endproclaim

\proclaim{ Lemma \lemCharneyR } {\rm[\refCharney; Prop.~3.3].}
Suppose that $X=A_1\cdot A_2\cdot\dots\cdot A_r$
is in left normal form
{\rm(}$A_i\in{]1,\Delta[}$, $i=1,\dots,r${\rm)},
and let $A_{r+1}$ be a simple element.
Then the decomposition $XA_{r+1}=A''_1\cdot\dots\cdot A''_{r+1}$is left weighted  where
the $A''_i$'s are defined recursively together with simple elements $A'_1,\dots,A'_r$ by
the conditions that
$A'_{r+1}=A_{r+1}$, $A'_i\cdot A''_{i+1}$ is the left normal form of
$A_i,A'_{i+1}$ for $i=r,\dots,1$, and $A''_1=A'_1$.
We have  $A''_i\ne\Delta$ for $i>1$ and $A''_i\ne1$ for $i\le r$ {\rm(}but
it is possible that $A''_1=\Delta$ or $A''_{r+1}=1${\rm)}.
\qed\endproclaim

\head\sectSSS. Super summit set of a product of two conjugates of atom powers
               in square-free homogeneous symmetric Garside groups
\endhead

In this section we prove Theorem \thSSS\ and Corollary  \corUSS.
Throughout this section $(G,\Cal P,\delta)$ is a square free symmetric homogeneous Garside structure with
set of atoms $\Cal A$.

\subhead\sectSSS.1. Preliminaries
\endsubhead

\proclaim{ Lemma \lemAx } {\rm[\refQPtwo; Lemma 3.1].}
Let $x\in\Cal A$ and $A\in\Cal P$. If $xA\preccurlyeq\delta$, {\rm(}resp. $Ax\preccurlyeq\delta${\rm)},
then there exists $x_1\in x^G\cap\Cal A$ such that $xA=Ax_1$ {\rm(}resp. $Ax=x_1A${\rm)}.
\endproclaim

\demo{ Proof } Immediately follows from the fact that the Garside structure is
symmetric and homogeneous.
\qed\enddemo

The following three results are proven in [\refQPtwo].

\proclaim{ Theorem \thQPone } {\rm[\refQPtwo; Theorem 1a].}
Let $X\sim x^k$ where $x\in\Cal A$, $k\ge1$. Then the left normal form of $X$ is
$\delta^{-n}\cdot A_n\cdot\dots\cdot A_1\cdot x_1^k\cdot B_1\cdot\dots\cdot B_n$ where
$n\ge0$, $x_1\in x^G\cap\Cal A$, and $A_i\delta^{i-1}B_i=\delta^i$ for $i=1,\dots,n$.
In particular, $\ell(X)=k+2n=k-2\inf X$.
\qed
\endproclaim

\proclaim{ Theorem \thBlock } {\rm(Blocking property [\refQPtwo; Corollary 7.2])}.
Let $X\sim x^k$ where $x\in\Cal A$, $X\not\in\Cal P$, $k\ge1$. 
Let $U\in G$ be such that $\inf XU = \inf X+\inf U$. Then $\iota(XU)=\iota(X)$.
\qed\endproclaim

\proclaim{ Lemma \lemAAP } {\rm[\refQPtwo; Lemma 7.5]}.
Let $A\in[1,\delta]$ and $P\in\Cal P$. Then $\delta\wedge(AP)=\delta\wedge(A^2P)$.
In particular, if $X\in G$ is such that $\inf AX=\inf X$, then $\iota(A^2X)=\iota(AX)$
and $\inf A^2X=\inf AX=\inf X$. \qed
\endproclaim

\smallskip\noindent
{\bf Remark \remAAP.} The conclusion of [\refQPtwo; Lemma 7.5] was erroneously stated in the form
$\iota(AP)=\iota(A^2P)$. This is wrong in general without the assumption $\delta\not\preccurlyeq AP$
as one can see in the example $G=\Br_4$ (with the Birman--Ko--Lee Garside structure,
thus $\delta=\sigma_3\sigma_2\sigma_1$),
$A=\sigma_2\sigma_1$, $P=\tau^2(A)$, and hence $\iota(AP)=\sigma_2$, $\iota(A^2P)=A$.
The statement and the proof of [\refQPtwo; Lemma 7.5] become correct
if one replaces all $\iota(\dots)$ by $\delta\wedge(\dots)$.
This mistake does not affect the usage of the lemma in the proof of the blocking property.

\proclaim{ Lemma \lemXt }
Let $x\in\Cal A$, $k\ge1$, $X\in(x^k)^G$, $s\in[1,\delta]$.
If $\ell(Xs)\le\ell(X)$ or $\ell(s^{-1}X)\le\ell(X)$, then $\ell(X^s)\le\ell(X)$.

\endproclaim

\demo{ Proof } If $\ell(Xs)\le\ell(X)$, then
$\ell(X^s)=\ell(s^{-1}Xs)\le\ell(s^{-1})+\ell(Xs)\le 1+\ell(Xs)\le1+\ell(X)$.
We have also $\ell(X^s)\equiv k\equiv\ell(X)\mod 2$ by
Theorem \thQPone. Hence $\ell(X^s)\le\ell(X)$.
The case $\ell(s^{-1}X)\le\ell(X)$ is similar.
\qed\enddemo


\proclaim{ Lemma \lemUXt }
Let $x\in\Cal A$, $k\ge1$, $X\in(x^k)^G$, $U\in G$, $s\in[1,\delta]$.
Suppose that
$$
       \sup UXs \le \sup UX = \sup U + \sup X.                    \eqno(\eqUXt)
$$
Then $\ell(X^s)\le\ell(X)$.
%
\endproclaim

\demo{ Proof }
The case $s\in\{1,\delta\}$ is trivial, so we assume that $s\in{]1,\delta[}$.
By Lemma \lemXt, it is enough to show that
$\sup Xs \le \sup X$. Suppose the contrary:
$$
     \sup Xs = \sup X + \sup s.                                \eqno(\eqUXtI)
$$
The inequality in (\eqUXt) can be rewritten as $\sup UXs < \sup UX + \sup s$.
By combining it with (\eqUXtI) and the equality in (\eqUXt), we obtain
$$
   \sup UXs < \sup UX + \sup s = \sup U + \sup X + \sup s = \sup U + \sup Xs.
$$
By Lemma \lemSup, this implies $\varphi(U)\phiRev(Xs)\preccurlyeq\delta$.
By Corollary \corCharney\ combined with (\eqUXtI) and Theorem \thQPone, we have $\phiRev(Xs)=\phiRev(X)$.
Hence $\varphi(U)\phiRev(X)\preccurlyeq\delta$
which contradicts  the equality in (\eqUXt).
\qed\enddemo


%

\subhead\sectSSS.2. Products of two atoms. Normal forms in Case (iii) of Theorem \thSSS
\endsubhead

\noindent
Recall that
$(G,\Cal P,\delta)$ is a square free symmetric homogeneous Garside structure with
set of atoms $\Cal A$.

\proclaim{ Proposition \propTwoAtoms }
Let $x$ and $y$ be two atoms such that $xy\preccurlyeq\delta$. Then there exist $m\ge 2$ and
pairwise distinct atoms $a_1,\dots,a_m$ {\rm(we assume that the indices are defined mod $m$)} such that:
\roster
\item"(i)" $x=a_1$, $y=a_2$, and $a_ia_{i+1}=xy$ for any $i$;
\item"(ii)" $a_{i+2} = a_i^{xy}$ for any $i$;
\item"(iii)" the product $a_i\cdot a_j$ is left weighted unless $j\equiv i+1\mod m$.
\endroster
\endproclaim

\demo{ Proof } We define $a_1,a_2,\dots$ recursively by $a_1=x$, $a_2=y$, $a_i a_{i+1}=a_{i-1} a_i$.
Then all $a_i$'s are atoms by Lemma \lemAx\ and (i) holds; (ii) follows from (i). Let us prove (iii).
Suppose that $a_i\cdot a_j$ is not left weighted, i.e., $a_ia_j\preccurlyeq\delta$.
Note that $a_i\vee a_j=xy$.
Since the Garside structure is symmetric, we have $a_i\prec a_ia_j$ and $a_j\prec a_i a_j$.
Hence $xy=a_i\vee a_j\preccurlyeq a_i a_j$. Since $\|xy\|=\|a_i a_j\|$, it follows that
$a_i a_j=xy=a_i a_{i+1}$ whence $a_j=a_{i+1}$. \qed
\enddemo

For $x,y\in\Cal A$, we set
$$
   \mu_{x,y}=\cases  0, &\text{if $x\cdot y$ is left weighted},\\
                     1, &\text{if $x=y$},\\
                     m, &\text{if $xy\preccurlyeq\delta$
                               and $m$ is as in Proposition \propTwoAtoms}.
                \endcases
$$

\noindent
{\bf Remark \remTwoAtoms.} It follows from Proposition \propTwoAtoms\ that the
submonoid of $G$ generated by any pair of atoms is either free or isomorphic to the positive monoid
of the dual Garside structure in an Artin-Tits group of type $I_2(m)$ (see [\refPicant; Proposition 1.2]).
It is interesting to study if the same is true for the subgroup of $G$ generated by a pair of atoms.
Note that the subgroup generated by a submonoid $M$ of a group is not necessarily isomorphic to the group of
fractions of $M$. For example, the submonoid $M$ of $\Br_3$ generated by $\sigma_1$ and $\sigma_2^{-1}$
is free whereas the subgroup generated by $M$ is the whole $\Br_3$ which is not a free group.

\proclaim{ Proposition \propLNF }
(a). Let $Z=x^ky^l$ where $k,l\ge 1$ and $x,y\in\Cal A$, $x\ne y$.
Then $Z\not\in\SSS(Z)$ if and only if one of the following conditions holds:
\roster
\item"(i)"   $\mu_{y,x}\ge 3$;
\item"(ii)"  $\mu_{x,y}=3$, $k=1$, and $l\ge 3$;
\item"(iii)" $\mu_{x,y}=3$, $l=1$, and $k\ge 3$.
\endroster
If $Z\in\SSS(Z)$, then the left normal form of $Z$ is
$$
   \cases
      x^k\cdot y^l                     &\text{if $\mu_{x,y}=\mu_{y,x}=0$},\\
      (xy)^k\cdot y^{l-k}              &\text{if $\mu_{x,y}=2$ and $k\le l$
                                              {\rm(the case $l\le k$ is similar)}},\\
      xy\cdot (x^y)^{k-1}\cdot y^{l-1} &\text{if $\mu_{x,y}\ge 3$.}
   \endcases
$$

\smallskip\noindent
(b). Let $Z=x^p y^l x^q$ where $p,q,l\ge 1$ and $x,y\in\Cal A$, $xy\ne yx$. Then
$Z\not\in\SSS(Z)$ if and only if one of the following conditions holds:
\roster
\item"(i)"  $\mu_{x,y}=3$, $p=l=1$, and $q\ge 2$;
\item"(ii)" $\mu_{y,x}=3$, $q=l=1$, and $p\ge 2$.
\endroster
If $Z\in\SSS(Z)$, then the left normal form of $Z$, is
$$
   \cases
      x^p\cdot y^l\cdot x^q                     &\text{if $\mu_{x,y}=\mu_{y,x}=0$},\\
      xy\cdot x_1^{p-1}\cdot y^{l-1}\cdot x^q   &\text{if either $\mu_{x,y}\ge 4$, or $\mu_{x,y}=3$ and $l\ge2$,}\\
      yx\cdot x_2^p\cdot y_2^{l-1}\cdot x^{q-1} &\text{if either $\mu_{y,x}\ge 4$, or $\mu_{y,x}=3$ and $l\ge2$,}\\
      (xy)^2\cdot y^{p-2}\cdot x^{q-1}   &\text{if $\mu_{x,y}=3$ and $l=1$,}\\
      (yx)^2\cdot y_2^{p-1}\cdot x^{q-2} &\text{if $\mu_{y,x}=3$ and $l=1$}\\
   \endcases
$$
where $x_1$, $x_2$, and $y_2$ are defined by $xy=yx_1$ and $yx=xy_2=y_2x_2$.
\endproclaim

\demo{ Proof } A straightforward computation using Proposition \propTwoAtoms. To see that
the listed elements $Z$ are in the super summit set, it is enough to check that in each case
$\frak s(Z)$ belongs to the same list and $\ell(\frak s(Z))=\ell(Z)$. Thus
$\ell(\frak s^m(Z))=\ell(Z)$ for any $m$ whence $Z\in\SSS(Z)$ by [\refGebGM].
\qed\enddemo


\subhead\sectSSS.3. Proof of Theorem \thSSS\ and Corollary \corUSS
\endsubhead
Recall that
$(G,\Cal P,\delta)$ is a square free symmetric homogeneous Garside structure with
set of atoms $\Cal A$.

For $x,y\in\Cal A$ and $k,l\ge1$, we set:
$$
\split
  &{\vec\Cal G}'_{p,q}(x^k,y^l) = \{XY\mid X\sim x^k,\,Y\sim y^l,\,\ell(X)=2p+k,\,\ell(Y)=2q+l,\\
  &\qquad\qquad\qquad\qquad\quad\;\ell_s(XY)=\ell(X)+\ell(Y)\},\\
  &{\vec\Cal G}''_{p,n}(x^k,y^l) =
       \{Z=x_1^pYx_1^{k-p}\mid Y\sim y^l,\, x_1\in x^G\cap\Cal A,\,\ell(Y)=2n+l,\\
  &\qquad\qquad\qquad\qquad\qquad\qquad\qquad\ell_s(Z)=k+\ell(Y)\},\\
  &{\vec\Cal G}'''_p(x^k,y^l) = \{Z=x_1^p y_1^l x_1^{k-p}\mid
   x_1\in x^G\cap\Cal A,\,   y_1\in y^G\cap\Cal A,\,   Z\in\SSS(Z)\}
\endsplit
$$
and $\Cal G(x^k,y^l) = \Cal G'(x^k,y^l) \cup \Cal G''(x^k,y^l) \cup \Cal G'''(x^k,y^l)$ where
$$
\split
  &{\vec\Cal G}'(\;)=\bigcup_{p,q\ge0}{\vec\Cal G}'_{p,q}(\;),\quad
   {\vec\Cal G}''(\;)=\bigcup_{0\le p\le k;\,n\ge0}{\vec\Cal G}''_{p,n}(\;),\quad
   {\vec\Cal G}'''(\;)=\bigcup_{0\le p\le k}{\vec\Cal G}'''_{p}(\;),
\\
   &{\Cal G}^*(x^k,y^l) = {\vec\Cal G}^*(x^k,y^l)\cup{\vec\Cal G}^*(y^l,x^k)
   \qquad\text{where ${}^*$ stands for ${}'$ or ${}''$ or ${}'''$}.
\endsplit
$$
It is clear that $Z\in\Cal G(x^k,y^l)$ implies $Z\in\SSS(Z)$.
In this notation, the conclusion of Theorem \thSSS\ reads as $\SSS(Z)\subset\Cal G(x^k,y^l)$.
Let us fix $k,l\ge1$ and $x,y\in\Cal A$.

\proclaim{ Lemma \lemG }
Let $Z\in{\vec\Cal G}'(x^k,y^l)$ 
and let $s$ be a simple element
such that $Z^s\in\SSS(Z)$. Then $Z^s\in\Cal G(x^k,y^l)$.
\endproclaim

\demo{ Proof } Let $Z=XY$, $X\sim x^k$, $Y\sim y^l$, $\ell(Z)=\ell(X)+\ell(Y)$.
Since $Z,Z^s\in\SSS(Z)$, we have $\ell(Z)=\ell(Z^s)$, hence
$\ell(X^s)+\ell(Y^s)\ge\ell(X^sY^s)=\ell(Z^s)=\ell(Z)$.
On the other hand, we have $\ell(X^s)\le\ell(s^{-1})+\ell(X)+\ell(s)=\ell(X)+2$
and, similarly, $\ell(Y^s)\le\ell(Y)+2$. We have also
$\ell(X^s)\equiv k\equiv\ell(X)$ and
$\ell(Y^s)\equiv l\equiv\ell(Y)\mod2$ by Theorem \thQPone. Hence
$$
    \ell(Z) \le \ell(X^s)+\ell(Y^s) \le \ell(Z)+4,
    \qquad
    \ell(X^s)+\ell(Y^s)\equiv \ell(Z)\mod2.
$$
Thus $\ell(X^s)+\ell(Y^s)$ may take only three values: $\ell(Z)$, $\ell(Z)+2$, and $\ell(Z)+4$.
We consider separately these three cases.

\smallskip
Case 1. $\ell(X^s)+\ell(Y^s) = \ell(Z)$. The result immediately follows.

\smallskip
Case 2. $\ell(X^s)+\ell(Y^s) = \ell(Z)+2$.
Then, for $(U,V)=(X,Y)$ or $(Y,X)$, we have $\ell(U^s)=\ell(U)$ and $\ell(V^s)=\ell(V)+2$,
hence $\inf U^s=\inf U$, $\sup U^s=\sup U$, $\inf V^s=\inf V-1$, $\sup V^s=\sup V+1$ and
we obtain
$$
   \inf X^s + \inf Y^s = \inf Z^s - 1
   \qquad\text{and}\qquad
   \sup X^s + \sup Y^s = \sup Z^s + 1.                               \eqno(\eqLemGii)
$$

\smallskip
Case 2.1. $\inf X^s=0$ or $\inf Y^s=0$. Without loss of generality we may assume that
$\inf X^s=0$, i.~e., $X^s=x_1^k$ where $x_1\in x^G\cap\Cal A$. In this case we have
$\ell(X^s)=\ell(X)$ and $\ell(Y^s)=\ell(Y)+2$.
Let $(A,B)=(\iota(Y^s),\,\varphi(Y^s))$.
Then, by Theorem \thQPone, we have $Y^s=A\delta^{-1}Y_1B$ with
$\ell(Y_1)=\ell(Y^s)-2=\ell(Y)$, $BA=\delta$, and hence, $Y^s=Y_1^B$.
By (\eqLemGii) combined with Lemma \lemSup, we have
$\delta\preccurlyeq\iotaRev(X^s)\iota(Y^s)$. Since $\iotaRev(X^s)=x_1$, we obtain
$\delta\preccurlyeq x_1A$. Since, moreover, $\|\delta\|\ge\|x_1\|+\|A\|$,
this yields $x_1A=\delta$. Since $BA=\delta$, we obtain $B=x_1$, hence
$Z^s = x_1^k Y^s=x_1^k Y_1^B=x_1^{k-1} Y_1 x_1$.
Since $Y_1\sim y^l$ and $\ell(Y_1)=\ell(Y)$, we conclude that $Z^s\in\Cal G(x^k,y^l)$.

\smallskip
Case 2.2.  $\inf X^s<0$ and $\inf Y^s<0$. Let
$(A,B)=(\phiRev(X^s),\iotaRev(X^s))$ and $(C,D)=(\iota(Y^s),\varphi(Y^s))$.
Then, by Theorem \thQPone, we have
$X^s=A\delta^{-1}X_1B$ and $Y^s=C\delta^{-1}Y_1D$ where $BA=DC=\delta$, $X_1\sim X$, $Y_1\sim Y$,
$\ell(X_1)=\ell(X^s)-2$, and $\ell(Y_1)=\ell(Y^s)-2$.
By (\eqLemGii) combined with Lemma \lemSup\ we have
$\iotaRev(X^s)\iota(Y^s)=E\delta$ for some $E\in[1,\delta]$.
Hence
$Z^s = A\delta^{-1}X_1BC\delta^{-1}Y_1D = A\delta^{-1}X_1 E Y_1D = \delta^{-1}\tilde A X_1 E Y_1 D$
where $\tilde A=\tau^{-1}(A)$. Since $\tilde AB=C\tau(D)=\delta$, we have
$\delta^2=\tilde ABC\tau(D)=\tilde AE\delta\tau(D)=\tilde AED\delta$ whence $\tilde AED=\delta$.

Case 2.2.1. $\ell(\tilde AX_1)\le\ell(X_1)$ or $\ell(Y_1D)\le\ell(Y_1)$. By symmetry, it
is enough to consider only the latter case. So, let $\ell(Y_1D)\le\ell(Y_1)$. Then,
by Lemma \lemXt, we have $\ell(Y_1^D)\le\ell(Y_1)$. Since
$Z^s = \delta^{-1}\tilde A X_1 ED Y_1^D= X_1^{ED} Y_1^D$ and
$\ell(X_1^{ED})+\ell(Y_1^D)\le (\ell(X_1)+2)+\ell(Y_1)=\ell(X^s)+(\ell(Y^s)-2)=\ell(Z^s)$,
we conclude that $Z^s\in\Cal G(x^k,y^l)$.

Case 2.2.2. $\ell(\tilde AX_1)=\ell(X_1)+1$ and $\ell(Y_1D)=\ell(Y_1)+1$.
Let us show that this is impossible.
Indeed, in this case we have $\sup \tilde AX_1=\sup\tilde A+\sup X_1=\sup X_1+1=\sup X^s$ and
similarly
$\sup Y_1D=\sup Y^s$. By (\eqLemGii), this yields
$$
   \sup \tilde AX_1+\sup Y_1D=\sup X^s+\sup Y^s=\sup Z^s+1 = \sup \tilde AX_1EY_1D.
$$
By Lemma \lemXtY, this implies that there exist $u,v\in[1,\delta]$ such that
$E=uv$, $\sup \tilde AX_1u=\sup \tilde AX_1$, and $\sup vY_1D=\sup Y_1D$.
Then, by Lemma \lemUXt, we have $\ell(X_2)\le\ell(X_1)$ and $\ell(Y_2)\le\ell(Y_1)$
where $X_2=u^{-1}X_1u$ and $Y_2=vY_1v^{-1}$.
Since
$$
   Z^s = \delta^{-1}\tilde AX_1uvY_1D=\delta^{-1}\tilde AuX_2Y_2vD=(X_2Y_2)^{vD},
$$
we obtain $\ell_s(Z)\le\ell(X_2Y_2)\le\ell(X_2)+\ell(Y_2)\le\ell(X_1)+\ell(Y_1)
=\ell(X^s)+\ell(Y^s)-4=\ell(Z^s)-2$. Contradiction.

\smallskip
Case 3. $\ell(X^s)+\ell(Y^s) = \ell(Z)+4$. Let us show that this case is impossible.
We have
$\ell(s^{-1}Xs)=\ell(s^{-1})+\ell(X)+\ell(s)$ and
$\ell(s^{-1}Ys)=\ell(s^{-1})+\ell(Y)+\ell(s)$, hence
$$
   \ell(s^{-1}X)=\ell(s^{-1})+\ell(X)
   \qquad\text{and}\qquad
   \ell(Ys)=\ell(Y)+\ell(s)                             \eqno(\eqLemGiii)
$$
whence
$\sup s^{-1}X=\sup s^{-1} + \sup X=\sup X$
and $\sup Ys=\sup Y + \sup s=\sup Y + 1$.
Thus
$\sup s^{-1}X + \sup Ys = \sup X +\sup Y+1>\sup X +\sup Y=\sup Z=\sup Z^s =\sup s^{-1}XYs$.
By Lemma \lemSup, this implies
$\varphi(s^{-1}X)\phiRev(Ys)\preccurlyeq\delta$.
We have $\varphi(s^{-1}X)=\varphi(X)$ by (\eqLemGiii) combined with Corollary \corCharney.
Similarly, $\phiRev(Ys)=\phiRev(Y)$.
Thus we obtain $\varphi(X)\phiRev(Y)\preccurlyeq\delta$ which contradicts the condition
$\ell(XY)=\ell(X)+\ell(Y)$.
\qed\enddemo

\proclaim{ Lemma \lemGp }
Let $Z\in{\vec\Cal G}''(x^k,y^l)$ and let 
$s$ be a simple element
such that $Z^s\in\SSS(Z)$. Then $Z^s\in\Cal G(x^k,y^l)$.
\endproclaim

\demo{ Proof }
Let $Z=x_1^pYx_1^q$ where $x_1\in x^G\cap\Cal A$, $Y\sim y^l$, $p+q=k$, $\ell(Z)=\ell(Y)+k$.
If $p=0$ or $q=0$, then Lemma \lemG\ applies. So, we assume that $p,q>0$.
Let us show that
$$
    \sup s^{-1}Z < \sup s^{-1} + \sup Z
    \qquad\text{or}\qquad
    \sup Zs < \sup Z + \sup s.                                   \eqno(\eqLemGpI)
$$
Indeed, suppose that the left inequality in (\eqLemGpI) does not hold, i.~e.,
$\sup s^{-1}Z = \sup s^{-1} + \sup Z=\sup Z$.
Then
$$
    \sup s^{-1}Z +\sup s = \sup Z+1 > \sup Z=\sup(s^{-1}Z\cdot s).
$$
Hence $\varphi(s^{-1}Z)s\preccurlyeq\delta$ by Lemma \lemSup.
Since $\varphi(s^{-1}Z)=\varphi(Z)$ by Corollary \corCharney,
this means that $\varphi(Z)s\preccurlyeq\delta$ which implies the right inequality in
(\eqLemGpI). Thus, (\eqLemGpI) holds.

By symmetry, without loss of generality we may assume
that the right inequality in (\eqLemGpI) holds.
Then $x_1s=\varphi(Z)s\preccurlyeq\delta$ by Lemma \lemSup.
Hence, by Lemma \lemAx, we have
$x_1s=sx_2$ where $x_2=x_1^s\in x^G\cap\Cal A$, and we obtain
$Z^s = x_2^p Y^s x_2^q$. If $\ell(Y^s)\le\ell(Y)$, then we are done. So, we suppose
that $\ell(Y^s)=\ell(Y)+2$. In this case we have also $\inf Y^s = \inf Y -1$.

Let us show that
$$
    \inf x_2^p Y^s > \inf x_2^p + \inf Y^s
    \qquad\text{or}\qquad
    \inf Y^s x_2^q > \inf Y^s + \inf x_2^q.                       \eqno(\eqLemGpII)
$$
Indeed, suppose that the right inequality in (\eqLemGpII) does not hold, i.~e.,
$\inf Y^s x_2^q = \inf Y^s + \inf x_2^q$,
hence $\inf x_2^p + \inf Y^s x_2^q = \inf x_2^p + \inf Y^s + \inf x_2^q =\inf Y^s <
\inf Y = \inf Z = \inf Z^s$.
Then we have $\delta\preccurlyeq\iotaRev(x_2^p)\iota(Y^s x_2^q)$ by Lemma \lemInf.
By Theorem \thBlock, we have $\iota(Y^s x_2^q)=\iota(Y^s)$.
Hence  $\delta\preccurlyeq\iotaRev(x_2^p)\iota(Y^s)$ which implies the left inequality in
(\eqLemGpII). Thus, (\eqLemGpII) holds.

By symmetry, without loss of generality we may assume
that the left inequality in (\eqLemGpII) holds.
The rest of the proof is almost the same as in Case 2.1 of Lemma \lemG. Namely,
let $(A,B)=(\iota(Y^s),\,\varphi(Y^s))$.
Then, by Theorem \thQPone, we have $Y^s=A\delta^{-1}Y_1B$ with
$\ell(Y_1)=\ell(Y^s)-2=\ell(Y)$, $BA=\delta$, and hence, $Y^s=Y_1^B$.
Then we have $\delta\preccurlyeq\iotaRev(x_2^p)\iota(Y^s)=x_2A$ by Lemma \lemInf\
combined with the left inequality in (\eqLemGpII).
Since $BA=\delta$, we obtain $B=x_2$, hence
$Z^s = x_2^p Y^s x_2^q = x_2^p Y_1^B x_2^q = x_2^{p-1} Y_1 x_2^{q+1}$.
Since $Y_1\sim y^l$ and $\ell(Y_1)=\ell(Y)$, we conclude that $Z^s\in\Cal G(x^k,y^l)$.
\qed\enddemo

\proclaim{ Lemma \lemGpp }
Let $Z\in{\vec\Cal G}'''(x^k,y^l)$ and let 
$s$ be a simple element
such that $Z^s\in\SSS(Z)$. Then $Z^s\in{\Cal G}'''(x^k,y^l)$.
\endproclaim

\demo{ Proof } We shall assume that $\|\delta\|\ge3$. In the case $\|\delta\|=2$,
the proof is the same but the notation should be slightly changed.

By the same arguments as in the proof of Lemma \lemGpp, we may
assume that the right inequality in (\eqLemGpI) holds.
By Proposition \propLNF, we have $\|\varphi(Z)\|=1$ or $2$.

\smallskip
Case 1. $\|\varphi(Z)\|=1$. It follows from 
Proposition \propLNF\ that, up to exchange of the roles of $x^k$ and $y^l$,
we may assume that $Z=x_1^p Y x_1^q$ where
$Y=y_1^l$, $x_1\in x^G\cap\Cal A$, $y_1\in y^G\cap\Cal A$, $p+q=k$, $q\ge1$, and $\varphi(Z)=x_1$.
The rest of the proof is the same as in Lemma \lemGp.

Note that the presentation of $Z$ in the form as in the definition of ${\Cal G}'''(x^k,y^l)$ is not necessarily
unique. For example, if $k=4$, $l=1$, and $Z=xyx^3$ where $xy=yz=zx$, $z\in\Cal A$, then
we work with $Z = x^1 y^1 x^3$, $\varphi(Z)=x$ when the right equality in (\eqLemGpI) holds,
but we work with $Z = y^4 z^1 y^0$, $\phiRev(Z)=y$ when the left equality in (\eqLemGpI) holds.

\smallskip
Case 2. $\|\varphi(Z)\|=2$. By Proposition \propLNF, we may assume that
$Z=x_0^p y_0^l x_0^q$, $p+q=k$, $x_0\in x^G\cap\Cal A$, $y_0\in y^G\cap\Cal A$,
and $\varphi(Z)=uv$ where $(u,v)$ is $(x_0,y_0)$ or $(y_0,x_0)$.
By the right inequality in (\eqLemGpI) combined with Lemma \lemSup, we have $\varphi(Z)s\preccurlyeq\delta$,
thus $uvs\preccurlyeq\delta$. Hence $vs\preccurlyeq\delta$ and $vs=sv_1$, $v_1=v^s\in\Cal A$ by Lemma \lemAx.
Then we have $usv_1=uvs\preccurlyeq\delta$ whence $us\preccurlyeq\delta$ and
$us=su_1$, $u_1=u^s\in\Cal A$. Thus $x_0^s=x_1$ and $y_0^s=y_1$ with $x_1,y_1\in\Cal A$, and
we obtain $Z^s = x_1^p y_1^l x_1^q \in{\Cal G}'''(x^k,y^l)$. \qed
\enddemo

\demo{ Proof of Theorem \thSSS }
As we already pointed out before Lemma \lemG, we need to prove that $\SSS(Z)\subset\Cal G(x^k,y^l)$.
We have $\SSS(Z)\cap\Cal G(x^k,y^l)\ne\varnothing$.
Indeed, if $Z\not\in\Cal P$, then $\SSS(Z)\cap{\vec\Cal G}''_k(x^k,y^l)\ne\varnothing$ by
[\refQPtwo; Theorem 1b] (in fact, only [\refQPtwo; Corollary 3.5] is needed here).
If $Z\in\Cal P$, then, again by [\refQPtwo; Theorem 1b], we have
$Z\sim Z_1=x_1^k y_1^l$ where $x_1\in x^G\cap\Cal A$, $y_1\in y^G\cap\Cal A$.
By Proposition \propLNF a, it follows that $Z_1\in\SSS(Z)$, and hence $Z_1\in{\Cal G}'''(x^k,y^l)$,
unless one of Cases (i)--(iii) of Proposition \propLNF\ occur. However, in each of these three
cases, a cyclic permutation of the word $x_1^k y_1^l$ yeilds an element $Z_2$ of
$\SSS(Z)$. Then we have $Z_2\in\SSS(Z)\cap{\Cal G}'''(x^k,y^l)$.

By the convexity theorem [\refEM; Corollary 4.2], any element of $\SSS(Z)$ can be obtained
from any other by successive conjugations by simple elements. Thus the result follows from
Lemmas \lemG\ -- \lemGpp.
\qed\enddemo

The following proposition shows that the cycling operator acts on the sets
${\vec\Cal G}'_{p,q}(x^k,y^l)$ and ${\vec\Cal G}''_{p}(x^k,y^l)$
in the most natural and expected way.

\proclaim{ Proposition \propSSSc } If $p>0$, then

\noindent
$\bold c\big( {\vec\Cal G}'_{p,q}(x^k,y^l)\big) \subset
                    {\vec\Cal G}'_{p-1,q+1}(x^k,y^l)$ and
$\bold c\big( {\vec\Cal G}''_{p,n}(x^k,y^l)\big) \subset
                    {\vec\Cal G}''_{p-1,n}(x^k,y^l)$.
Note that

\noindent
${\vec\Cal G}'_{0,n}(x^k,y^l) =  {\vec\Cal G}''_{k,n}(x^k,y^l)$ and
${\vec\Cal G}''_{0,n}(x^k,y^l) =  {\vec\Cal G}'_{n,0}(y^l,x^k)$.
\endproclaim

\demo{ Proof }
The first inclusion follows from Corollary \corSSS. Let us prove the second one.
Let $Z$ be as in the definition of ${\vec\Cal G}'_{p,n}(x^k,y^l)$. We may suppose that
the left normal form of $Z$ is as in (\eqSSSii).
We see from (\eqSSSii) that $\iota(Z)=\iota(x_1^p Y)=\tilde C_1=\tau^n(C_1)$.
By Lemma \lemAAP, we have $\iota(x_1^pY)=\iota(x_1 Y)$.
Hence $\tilde C_1=x_1s=sx_2$ where $s\preccurlyeq\iota(Y)$ and $x_2\in x^G\cap\Cal A$.
Thus $Z = x_1^psY' x_1^{k-p} = sx_2^pY' x_1^{k-p} = \tilde C_1x_2^{p-1}Y' x_1^{k-p}$
and $\bold c(Z)=x_2^{p-1}Y' x_1^{k-p}\tilde C_1=x_2^{p-1}Y' x_1^{k-p+1}s=x_2^{p-1}Y's x_2^{k-p+1}
\in{\vec\Cal G}'_{p-1,n}(x^k,y^l)$.
\qed\enddemo

Corollary \corUSS\  follows from Proposition \propSSSc.



\head\sectQPfour. Homogeneous symmetric Garside groups with $\|\delta\|=3$
\endhead

In this section we assume that $(G,\Cal P,\delta)$ is a
square free homogeneous symmetric Garside structure with set of atoms $\Cal A$
and we assume that $\|\delta\|=3$.

If $\delta^p\cdot A_1\cdot\dots\cdot A_n$ is the left normal form of $X$, then we
denote:
$$
   \ell_1(X) = \operatorname{Card}\{i\mid\; \|A_i\|=1\},
   \qquad
   \ell_2(X) = \operatorname{Card}\{i\mid\; \|A_i\|=2\}.                 \eqno(\eqDefCount)
$$

\proclaim{ Lemma \lemCount } Let $X\in G$.
Then 
$$
     \ell_1(X) = \inf X +2\sup X - e(X)\qquad\text{and}\qquad
     \ell_2(X) = -2\inf X -\sup X + e(X).
$$
\endproclaim


\demo{ Proof } Follows from $n_1+n_2=\ell(X)$ and
$n_1 + 2n_2 = e(X)-3\inf X$, $n_i=\ell_i(X)$. \qed
\enddemo

\proclaim{ Lemma \lemC }
Let $Y=\delta^p\cdot A_1\cdot\dots\cdot A_n$ be in left normal form, $n\ge3$.
Suppose that $\inf_s Y > p$.

(a). If $\iota(\bold c(Y))=\tau^{-p}(A_2)$, then  $\inf\bold c(Y) > p$.

(b). If $(\|A_2\|,\dots,\|A_n\|)\ne(1,\dots,1)$, then $\inf\bold c(Y) > p$.
\endproclaim

\demo{ Proof } (a). If $\iota(\bold c(Y))=\tilde A_2$, then
$\bold c^2(Y)=\delta^p A_3\dots A_n\,\tilde A_1\tilde A_2$ where $\tilde A_j=\tau^{-p}(A_j)$.
Since $\inf_s Y > p$, it follows from [\refBKLtwo] that $\inf\bold c^2(Y) > p$.
 Hence $\delta\preccurlyeq A_3\dots A_n\,\tilde A_1\tilde A_2$.
Then, by Lemma \lemInf, we have $\delta\preccurlyeq\iotaRev(A_3\dots A_n) \tilde A_1$,
hence $\delta\preccurlyeq A_2\dots A_n\, \tilde A_1$ which means that
$\inf\bold c(Y)>p$.

\smallskip
(b).  Suppose that $(\|A_2\|,\dots,\|A_n\|)\ne(1,\dots,1)$.
Let  $i\ge2$ be such that $\|A_i\|=2$.
Suppose that $\inf\bold c(Y) = p$. Then, by Lemma \lemCharneyR,
the left normal form of $\bold c(Y)$ starts with
$\delta^p\cdot A_2\cdot\dots\cdot A_i$. Hence $\inf\bold c(Y)>p$ by (a). Contradiction
\qed\enddemo

\proclaim{ Lemma \lemD }
Let $Y=\delta^p\cdot A_1\cdot\dots\cdot A_n$ be in left normal form, $n\ge3$.
Suppose that $\sup_s Y < p+n$.

(a). If $\varphi(\bold d(Y))=A_{n-1}$, then $\sup\bold d(Y) < p+n$.

(b). If $(\|A_1\|,\dots,\|A_{n-1}\|)\ne(2,\dots,2)$, then $\sup\bold d(Y) < p+n$.
\endproclaim


\demo{ Proof } Apply Lemma \lemC\ to $Y^{-1}$. \qed\enddemo

\proclaim{ Lemma \lemCDi }

(a). Let $\inf Y < \inf\bold c(Y)$ and $\sup\bold c(Y) = \sup Y$. Then $\ell_2(Y)\ge 2$.

(b). Let $\inf Y = \inf\bold d(Y)$ and $\sup\bold d(Y) < \sup Y$. Then $\ell_1(Y)\ge 2$.
\endproclaim

\demo{ Proof } (a).
Let $A=\iota(Y)$, $Y=AY_1$, and $B=\iotaRev(Y_1)$. The condition $\inf Y < \inf\bold c(Y)=\inf Y_1A$ combined with
Lemma \lemInf\ implies $\delta\preccurlyeq BA$. The condition $\sup\bold c(Y) = \sup Y$
implies $\delta\ne BA$. Hence $\|BA\|>\|\delta\|=3$ whence
$\|B\|=\|A\|=2$.

\smallskip
(b).  Apply (a) to $Y^{-1}$.
\qed\enddemo

\proclaim{ Lemma \lemCDii } Let $\ell(Y)\ge 3$ {\rm(note that this is so when $e(Y)\ge2$ and $\inf Y<0$)}.

(a). If $\inf Y < \inf_s Y$ and $\sup_s Y = \sup Y$, then $\inf Y<\inf\bold c(Y)$.

(b). If $\inf Y = \inf_s Y$ and $\sup_s Y < \sup Y$, then $\sup\bold d(Y)<\sup Y$.

(c). If $Y\not\in\SSS(Y)$, then $\inf Y<\inf\bold c(Y)$ or $\sup\bold d(Y) < \sup Y$.
\endproclaim

\demo{ Proof }
(a). If $\inf Y < \inf_s Y$, then $\inf Y = \inf X < \inf\bold c(X)$ where $X=\bold c^m(Y)$ for some $m\ge0$
(see [\refBKLtwo]).
If, moreover, $\sup_s Y = \sup Y$, then $\ell_2(X)\ge 2$ by Lemma \lemCDi.
We have $\ell_2(X)=\ell_2(Y)$ by Lemma \lemCount, thus $\ell_2(Y)\ge2$, and
the result follows from Lemma \lemC b.

\smallskip
(b). Apply (a) to $Y^{-1}$.

\smallskip
(c). If $\inf Y = \inf_s Y$ or $\sup_s Y = \sup Y$, then the result follows from (a), (b).
Otherwise it follows from Lemmas \lemC b,  \lemD b because $\ell_2(Y)>1$ or $\ell_1(Y)>1$.
\qed\enddemo

\proclaim{ Lemma \lemUxVi } Let $Y\in a^G b^G$ where $a,b\in\Cal A$.
Suppose that $\inf_s Y<0$ and $\inf Y=\inf_s Y$ {\rm(}i.~e., $Y$ is in its summit set{\rm)}.
Then there exist $U,V\in G$ such that, up to exchange of $a$ and $b$, we have
$Y=UyV$ with $y\in a^G\cap\Cal A$, $UV\sim b$ and the following conditions hold:
$\ell(U)\ge1$, $\ell(V)\ge 1$, the product
$\varphi(U)\cdot y\cdot\iota(V)$ is left weighted, and hence
$\ell(Y) = \ell(U)+1+\ell(V)$.
\endproclaim

\demo{ Proof }
Induction on $\sup Y-\sup_s Y$.
If $\sup Y-\sup_s Y=0$, then $Y\in\SSS(Y)$, and the result
follows from Corollary \corSSS.
Indeed, if $Y=z^P y^Q$ with $\ell(Y)=2+2\ell(P)+2\ell(Q)$ and $\ell(Q)\ge1$, then we set
$U=z^PQ^ {-1}$ and $V=Q$; if $Y=y^Pz$ with $\ell(Y)=2+2\ell(P)$, then we set $U=P^{-1}$ and $V=Pz$.

Suppose that $\sup Y-\sup_s Y>0$. Then $\sup\bold d(Y)=\sup Y-1$ by Lemma \lemCDii{b}.
So, by the induction hypothesis, we assume that
$\bold d(Y)=U'y'V'$ with the required properties. Without loss of generality we may assume also
that $\inf V'=0$.

Let $\delta^p\cdot A_1\cdot\dots\cdot A_n$ be the left normal form of $Y$. Then the left
normal form $B_1\cdot\dots\cdot B_{n-1}$ of $\delta^{-p}\bold d(Y)$
is obtained from $\tau^p(A_n)\cdot A_1\cdot\dots\cdot A_{n-1}$
by the procedure described in Lemma \lemCharney. It follows that for some $i\ge1$, we have
$(\|A_n\|,\|A_1\|,\dots,\|A_{i-1}\|,\|A_i\|)=(1,2,\dots,2,1)$,
$(\|B_1\|,\dots,\|B_i\|)=(2,\dots,2)$, and $A_\nu=B_\nu$ for $\nu>i$; see Figure \figD.
Hence we have $U'=\delta^p B_1\dots B_{j-1}$, $y'=B_j$, $V'=B_{j+1}\dots B_{n-1}$
for some $j$ in the range $i<j<n-1$ and we obtain the desired decomposition
$Y=UyV$ by setting
$U=A_n^{-1}U'=\delta^p A_1\dots A_{j-1}$, $y=y'=A_j$, $V=V'A_n=A_{j+1}\dots A_n$.
\qed\enddemo

\midinsert
\centerline{
  $\tau^p(A_n)$ \hskip 1mm
  $A_1$ \hskip 9mm
  $\dots$ \hskip 9mm
  $A_{i-1}$ \hskip 1mm
  $A_i$ \hskip 15mm
  $A_i$ \hskip 1mm
  $A_{i+1}$ \hskip 2mm
  $\dots$ \hskip 2mm
  $A_n$\hskip 1mm
  $\tau^{-p}(A_1)$
}
\vskip 7mm
\centerline{\hskip 6mm
  $B_1$ \hskip 10mm
  $\dots$ \hskip 10mm
  $B_{i-1}$ \hskip 3mm
  $B_i$ \hskip 21mm
  $\delta$ \hskip 6mm
  $B_i$ 
  $B_{i+1}$ \hskip 1.5mm
  $\dots$ \hskip 1.5mm
  $B_n$
}
\vskip-11mm
\centerline{\epsfxsize=110mm\epsfbox{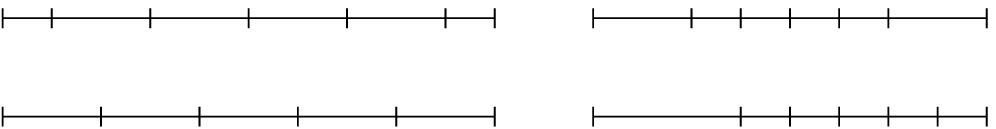}}
\vskip -2mm
\botcaption{\hskip 5mm Figure {\figD}  \hskip 40mm Figure \figC}
\endcaption
\endinsert

\proclaim{ Lemma \lemUxVs } Let $Y\in a^G b^G$ where $a,b\in\Cal A$.
Suppose that $\sup_s Y>1$, $\sup Y=\sup_s Y$ {\rm(}i.~e., $Y^{-1}$ is in its summit set{\rm)},
and $\|\varphi(Y)\|=1$.
Then there exist $U,V\in G$ such that, up to exchange of $a$ and $b$, we have
$Y=UyV$ with $y\in a^G\cap\Cal A$, $UV\sim b$ and the following conditions hold:
\roster
\item "(i)"   $\ell(V)\ge1$;
\item "(ii)"  $\ell(yV)=1+\ell(V)$;
\item "(iii)" if $\ell(U)>0$, then the product $\varphi(U)\cdot\iota(yV)$ is left weighted;
\item "(iv)"  if $\ell(U)>0$, then $\ell_2\big(\varphi(U)yV\big)\ge1$.
\endroster
Note that (ii) and (iii) imply $\ell(Y) = \ell(U)+1+\ell(V)$.
\endproclaim

\demo{Proof}
Induction on $\inf_s Y-\inf Y$. If $\inf_s Y-\inf Y=0$, then $Y\in\SSS(Y)$, and
the result follows from Corollary \corSSS.
Indeed, if $Y=y^P z^Q$ with $\ell(Y)=2+2\ell(P)+2\ell(Q)$, then we set
$U=P^{-1}$ and $V = P z^Q$.

Suppose that $\inf_s Y-\inf Y>0$. Then $\inf\bold c(Y)=\inf Y+1$ by Lemma \lemCDii{a}.
Let $\delta^p\cdot A_1\cdot\dots\cdot A_n$ be the left normal form of $Y$. We set $\tilde A_1=\tau^{-1}(A_1)$.
Then the left normal form $\delta\cdot B_2\cdot\dots\cdot B_n$ of $\delta^{-p}\bold c(Y)$
is obtained from $(A_2\cdot\dots\cdot A_n)\tilde A_1$ by the procedure described in Lemma \lemCharneyR:
$$
\split
  (A_2\cdot\dots\cdot A_n)\tilde A_1
  &=(A_2\cdot\dots\cdot A_{n-1})(C_n\cdot B_n)=\dots
  \\
  &=(A_2\cdot\dots\cdot A_{i-1}\cdot A_i)(C_{i+1}\cdot B_{i+1}\cdot\dots\cdot B_n)\\
  &=(A_2\cdot\dots\cdot A_{i-1})(\;\delta\;\cdot\; B_i\;\cdot\; B_{i+1}\cdot\dots\cdot B_n)=\dots\\
  &=(\delta\cdot B_2\cdot\dots\cdot B_{i-1}\cdot B_i\cdot\dots\cdot B_n)
\endsplit
$$
where $2\le i\le n$,
all the products in the parentheses are left weighted, and $B_\nu=\tau(A_\nu)$ for $\nu=2,\dots,i-1$.
It follows that
$(\|A_i\|,\|A_{i+1}\|,\dots,\|A_n\|,\|A_1\|)=(2,1,\dots,1,2)$ and
$(\|B_i\|,\dots,\|B_n\|)=(1,\dots,1)$; see Figure \figC.
Note that the condition $\|\varphi(Y)\|=1$ reads as $\|A_n\|=1$. Since $\|A_i\|=2$, this yields $i<n$.

Since $\varphi(\bold c(Y))=B_n$ and $\|B_n\|=1$, we may assume that the induction hypothesis holds,
so, we have a decomposition $\bold c(Y)=U'y'V'$ with the required properties. Without loss of generality we may assume also
that $\inf V'=0$.
We shall refer to Conditions (i)--(iv) applied to the decomposition $\bold c(Y)=U' y' V'$ by writing (i)$'$--(iv)$'$.
Condition (iii)$'$ means that
$U'=\delta^{p+1} B_2\dots B_{j-1}$ and $y'V'=B_j\dots B_n$ for some $j\ge2$.
Condition (iv)$'$ combined with $\|B_i\|=\dots\|B_n\|=1$ implies $j\le i$.

Let $U=\delta^p A_1\dots A_{j-1}$, $y=\tau^{-1}(y')$, and $V=y^{-1}A_j\dots A_n$.
First, let us show that $y\preccurlyeq A_j$. Indeed, if $j<i$, then $y=\tau^{-1}(y')\preccurlyeq\tau^{-1}(B_j)=A_j$.
If $j=i$, then $yC_{i+1}\preccurlyeq y\delta=\delta y'=\delta B_i=A_iC_{i+1}$ whence $y\preccurlyeq A_i=A_j$.
Thus
$$
      A_j=ys,\qquad V=s\cdot(A_{j+1}\cdot\dots\cdot A_n),\qquad s\in{[1,\delta[}.                    \eqno(\eqLemUxVs)
$$
We have $y\sim y'\sim a$ and $UV=\tilde A_1U'V'\tilde A_1^{-1}\sim U'V'\sim b$.
Let us show that the decomposition $Y=UyV$ satisfies (i)--(iv).
Indeed, $i<n$ implies (i), $\|A_i\|=2$ implies (iv),
and the fact that $A_1\cdot\dots\cdot A_n$ is left weighted implies (iii).
So, it remains to check that (ii) holds.
By (\eqLemUxVs) we have $\ell(V)\le\ell(yV)\le\ell(V)+1$, thus
it is enough to exclude the case $\ell(V)=\ell(yV)$, that is $\ell(V)=n-j+1$.

Suppose that $\ell(V)=n-j+1$. The product of $n-j$ factors in the parentheses in (\eqLemUxVs) is left weighted,
hence $A_n\succcurlyeq\varphi(V)$ by Lemma \lemCharney.
Since $A_n=\varphi(Y)$, we have $\|A_n\|=1$ by the hypothesis of the lemma. Thus the condition
$A_n\succcurlyeq\varphi(V)$ implies $A_n=\varphi(V)$.
We have
$$
\xalignat 2
   \sup V\tilde A_1 &= \sup V' +1  &&\;\text{because $V\tilde A_1 = \delta V'$}\\
                        &= \sup y'V'   &&\;\text{because $\ell(y'V')=\ell(V')+1$ by (ii)$'$}\\
                        &= n-j+1       &&\;\text{because $y'V'=B_j\dots B_n$}
\endxalignat
$$
hence $\sup V\tilde A_1=\sup V$ which implies $A_n\tilde A_1=\varphi(V)\tilde A_1\preccurlyeq\delta$
by Lemma \lemSup.
Hence $\sup\bold c(Y)<\sup Y$ which is impossible because $\sup Y=\sup_s Y$.
\qed\enddemo


\proclaim{ Lemma \lemForCaseII }
Let $V\in G$ and $x,y\in\Cal A$ be such that:
\roster
\item"(i)" $\ell(yV)=1+\ell(V)\ge 2$;
\item"(ii)" $\inf yVx=\inf yV$;
\item"(iii)" $\sup yVx=\sup yV$.
\endroster
Let $t=\varphi(yVx)$ and $yVx=Wt$.
Then $y\preccurlyeq\phiRev(W)$.
\endproclaim

\demo{ Proof }
Without loss of generality we may assume that $\inf yVx=\inf V=0$.
Then we have $\ell(U)=\sup U$ for elements $U$ of $G$ considered in this proof. Let $r=\ell(V)$.
The fact that $t=\varphi(Wt)$ implies $\ell(W)=\ell(Wt)-1$, hence
$$
    \ell(W)=\ell(yVx)-1=\ell(yV)-1=\ell(V)=r.                                              \eqno(\eqForCaseII)
$$
Let $A_1\cdot\dots\cdot A_r$ and $B_0\cdot\dots\cdot B_r$, $r\ge1$, be the left normal form of $V$ and of $yV$
respectively. By (ii) and (iii) we have $\delta\succ B_rx$.
Since $B_rx\succcurlyeq t$, we may write $B_rx=st$ with $s\in{[1,\delta[}$.
It follows from Lemma \lemCharneyR\ that $B_{r-1}s\cdot t$ is the left normal form of
$B_{r-1}\cdot B_rx$, in particular,
$$
    B_{r-1}s\preccurlyeq\delta  
                                                                                           \eqno(\eqForCaseIIa)
$$
Let $i$ be the minimal non-negative integer such that $A_j=B_j$ for all $j>i$.

%

\medskip
Case 1. $i=r$. Then we have $\|B_0\|=\dots=\|B_{r-1}\|=2$ by Lemma \lemCharney.
Hence the left normal form of $yVx$ is $B_0\cdot\dots\cdot B_{r-1}\cdot B_rx$.
Therefore the right normal form of $W$ is $B_0\cdot\dots\cdot B_{r-1}$,
and we obtain $y\preccurlyeq B_0 = \phiRev(W)$.

\medskip
Case 2. $i=r-1$ and $s=1$.
Then $t=A_rx=B_rx$ and $W=y\cdot A_1\cdot\dots\cdot A_{r-1}$, hence
$y=\phiRev(W)$ by (\eqForCaseII).

\medskip
Case 3. $i=r-1$ and $s\ne 1$. Then we have $\|B_0\|=\dots=\|B_{r-2}\|=2$ by Lemma \lemCharney.
Hence the left normal form of $yVx$ is $B_0\cdot\dots\cdot B_{r-2}\cdot B_{r-1}s\cdot t$
and the left normal form of $W$ is $B_0\cdot\dots\cdot B_{r-2}\cdot B_{r-1}s$.
The right normal form of $W$ coincides with
the left normal form because the letter length of each canonical factor is 2.
Hence $y\preccurlyeq B_0 = \phiRev(W)$.

\medskip
Case 4. $i\le r-2$.
Then $B_r=A_r$, $B_{r-1}=A_{r-1}$, and $W=y\cdot A_1\cdot\dots\cdot A_{r-2}\cdot B_{r-1}s$.
By (\eqForCaseIIa), this is a decomposition of $W$ into a product of $r$
simple elements. Hence $y=\phiRev(W)$ by (\eqForCaseII).
\qed\enddemo



\head\sectProof. Proof of Theorem \thMain
\endhead

Let the hypothesis of Theorem \thMain\ hold.
For a permutation $(\lambda,\mu,\nu)$ of $(1,2,3)$ and an integer
$n$, we set
$$
\split
  &\Cal Q_{n,p}^{(\lambda)} = \{ (x,Y)
  \mid  x Y\sim X,\;
  x\in a_\lambda^G\cap\Cal A,\; Y\in a_\mu^G a_\nu^G,\;
  \ell(Y) \le n,\; \inf Y \ge p\},
\\
  &\Cal Q_{n,p}= \Cal Q_{n,p}^{(1)}\cup \Cal Q_{n,p}^{(2)}\cup \Cal Q_{n,p}^{(3)},
\qquad \Cal Q_n={\bigcup}_p\Cal Q_{n,p}.
\qquad\text{and}\qquad\Cal Q={\bigcup}_n\Cal Q_{n}.
\endsplit
$$

Till the end of the section 
$(x,y,z)$ will always denote some permutation of $(b_1,b_2,b_3)$
with $b_i\in a_i^G\cap\Cal A$, and $x_1,x_2,\dots$ (resp. $y_1,y_2,\dots$ or $z_1,z_2,\dots$)
will stand for some atoms which are conjugate to $x$ (resp. to $y$ or to $z$).
All these new atoms will be obtained from $x$, $y$, $z$ by applying Lemma \lemAx.


\proclaim{ Lemma \lemCaseI }
Let $(x,Y)\in\Cal Q_{n,p}$, $p<0$.
Suppose that $\inf x Y > p$ or $\inf Y x > p$.
Then $\Cal Q_{n-1}\ne\varnothing$.
\endproclaim

\demo{ Proof } By symmetry, it is enough to consider the case when $\inf x Y > \inf Y$.
Let $A=\iota(Y)$. Then $\delta\preccurlyeq xA$ by Lemma \lemInf.
Since $\|x\|=1$ and $\|A\|\le2$, this means
$$
    xA = Ax_1 = \delta.                                            \eqno(\eqCaseI)
$$

Case 1. $Y\in\SSS(Y)$.
Then, by Corollary \corSSS, we have $Y=AUyV$ with $\ell(Y)=\ell(U)+\ell(V)+2$ and $AUV\sim z$.
Hence, for $Z=VxAU=V\delta U$, we have
$yZ=yVxAU \sim xAUyV=xY\sim X$ and
$Z = VxAU \sim xAUV \in x(z^G)$. Since
$\ell(Z)=\ell(V\delta U)\le\ell(V)+\ell(U)=\ell(Y)-2\le n-2$,
we obtain $(y,Z)\in\Cal Q_{n-2}$.

\medskip
Case 2. $Y\not\in\SSS(Y)$.
By Lemma \lemCDii c, $\inf Y < \inf\bold c(Y)$ or $\sup\bold d(Y) < \sup Y$.
If $\inf Y < \inf\bold c(Y)$, then $(xY)^A=x_1\bold c(Y)$ by (\eqCaseI), whence
$(x_1,\bold c(Y))\in\Cal Q_{n-1}$.

Suppose that $\sup\bold d(Y) < \sup Y$. Let $B=\varphi(Y)$, $Y=Y_1B$.
Then $\bold d(Y)=BY_1$ and $\ell(Y_1)=\ell(Y)-1$.
Let  $C=\phiRev(Y_1)$, $Y_1=CY_2$. Then $\ell(Y_2)=\ell(Y_1)-1=\ell(Y)-2$.
Since $\sup(BY_1) =\sup\bold d(Y) < \sup Y = \sup B + \sup Y_1$,
we obtain $BC\preccurlyeq\delta$ by Lemma \lemSup.
We have
$C=\phiRev(Y_1)\preccurlyeq\iota(Y_1)=\iota(Y)=A$ whence $xC\preccurlyeq xA=\delta$ by
(\eqCaseI). Hence $xC=Cx_2$ and we obtain
$(x Y)^C=x_2 Y^C$ with $\ell(Y^C)=\ell(Y_2BC)\le\ell(Y_2)+\ell(BC)=\ell(Y_2)+1=\ell(Y)-1$,
thus $(x_2,Y^C)\in\Cal Q_{n-1}$.
\qed\enddemo


\proclaim{ Lemma \lemCaseII }
Let $(x,Y)\in\Cal Q_{n,p}$ and $p<0$.
Suppose that $\sup x Y \le \sup Y$ or $\sup Y x \le \sup Y$.
Then either $x Y\in\SSS(X)$, or $Y x\in\SSS(X)$, or $\Cal Q_{n-1}\ne\varnothing$.
\endproclaim

\demo{ Proof }
By symmetry, it is enough to consider only the case $\sup Yx \le \sup Y$.
Then $Ax=x_1A\preccurlyeq\delta$ with $x_1\in x^G\cap\Cal A$ and $A=\varphi(Y)$, $Y=Y_1A$.
By Lemma \lemCaseI\ we may assume that
$$
     \inf xY = \inf Yx = \inf Y.                                 \eqno(\eqCaseIIinf)
$$
Let $B=\iotaRev(Yx)$.
Since the simple element $Ax$ divides $\delta^{-p}Yx$ from the right
but $\delta$ does not due to (\eqCaseIIinf), we conclude that
$B\succcurlyeq Ax$. Since $\|Ax\|=2$, this means that $B=Ax$.
Then $\cRev(Yx)=BY_1$.
If $B\cdot\iota(Y_1)$ is not left weighted, then $\inf BY_1 > \inf Y_1=p$
and the result follows from Lemma \lemCaseI\ applied to $(x_1,\bold d(Y))$ because
$x_1 \bold d(Y)=x_1AY_1 = BY_1$.
So, we assume that $B\cdot\iota(Y_1)$ is left weighted whence $\inf BY_1 = \inf Y_1$ which means that
$\inf\cRev(Yx) = \inf Yx$. By Lemma \lemC b this implies that either
$$
        \inf Yx = {\inf}_s Yx                                \eqno(\eqCaseII)
$$
or $\ell_2(Y)=0$.

\medskip
Case 1. $\ell_2(Y)=0$. Let $C=\iota(Y)$, $Y=CY_2A$. If $A\cdot C$ is left weighted,
then $Y$ is rigid, hence $Y\in\SSS(Y)$ which contradicts [\refQPtwo; Corollary 3]. Hence
$AC\preccurlyeq\delta$ and we obtain $(x_1,\bold d(Y))\in\Cal Q_{n-1}$ because
$x_1\bold d(Y) = x_1\bold d(Y_1A)=x_1AY_1=AxY_1\sim xY_1A=xY\sim X$ and
$\ell(\bold d(Y)) = \ell(\bold d(CY_2A))=\ell(ACY_2)\le\ell(AC)+\ell(Y_2)=1+\ell(Y_2)=\ell(Y)-1$.

\medskip
Case 2. $\ell_2(Y)>0$, thus (\eqCaseII) holds.
If $\sup_s Yx = \sup Yx$, then $Yx\in\SSS(X)$ and we are done. So, we assume that
${\sup}_s Yx < \sup Yx$ which implies by Lemma \lemCDii{b}
$$
       \sup\bold d(Yx) < \sup Yx.                                 \eqno(\eqCaseIIa)
$$

\medskip
Case 2.1. $\sup_s Y=\sup Y$. 
Suppose that $\sup_s Y\le 1$. Then $\inf_s Y=0$ and $\sup_s Y=1$ by [\refQPtwo; Theorem 1b]
(or by Corollary \corSSS).
By Lemma \lemCount, this yields
$p = \inf Y = e(Y) - 2\sup Y + \ell_1(Y) = 2 - 2\times 1 + \ell_1(Y) \ge 0$ which contradicts the hypothesis $p<0$.
Thus $\sup_s Y>1$. 
Recall also that $\varphi(Y)=A$ and $Ax\prec\delta$ whence $\|A\|=1$.

So, we may use
Lemma \lemUxVs. Hence $Y=UyV$ where $UV\sim z$ and Conditions (i)--(iv) of Lemma \lemUxVs\ hold.
Condition (iii) implies $\varphi(yV)=\varphi(Y)=A$.
Condition (iv) implies that the left normal form of $Vx$ coincides with the tail of the left normal form of
$Yx$, in particular, $\varphi(Yx)=\varphi(yVx)$; we denote this element by $t$ and we set $yVx=Wt$ as in Lemma \lemForCaseII.
Then we have $y\preccurlyeq\phiRev(W)$ by Lemma \lemForCaseII\
and we set $\phiRev(W)=ys=sy_1$, $W=ysW_1$with $s\in{[1,\delta[}$.

We are going to prove that $(y_1,Z)\in\Cal Q_{n-1}$ for $Z=W_1tUs$. We evidently have:
$$
\split
y_1Z&=y_1W_1tUs\sim sy_1W_1tU=ysW_1tU=WtU=yVxU\sim xUyV=xY\sim X,\\
Z&=W_1tUs\sim sW_1tU = VxU \sim xUV\in x^Gz^G.
\endsplit
$$
So, it remains to show that $\ell(Z)<n$.
We have $\bold d(Yx)=\bold d(UWt)=tU W$ and $\sup(\bold d(Yx))<\sup(Yx)=\sup(Y)$ by (\eqCaseIIa),
thus
$$
             \sup tU W < \sup Y.                    \eqno(\eqCaseIIb)
$$
If $\sup tU + \sup W < \sup Y$, then
$\ell(Z) = \ell(W_1tUs)\le\ell(W_1)+\ell(tU)+1=\ell(W)+\ell(tU)<\ell(Y)=n$ and we are done.
%
So, we assume that $\sup tU + \sup W \ge \sup Y$. Since
$\sup tU + \sup W\le 1+\sup U + \sup W=\sup U + \sup Wt =\sup U+\sup yVx = \sup U+\sup yV = \sup Y$,
it follows that $\sup tU + \sup W = \sup Y$.
Then (\eqCaseIIb) combined with Lemma \lemSup\
yields $\varphi(tU)\phiRev(W)\preccurlyeq\delta$ whence
$Bs\preccurlyeq Bsy_1=B\phiRev(W)\preccurlyeq\delta$ where $B=\varphi(tU)$. Thus, by setting $tU=U_1B$,
we obtain
$\ell(Z) = \ell(W_1tUs)\le\ell(W_1U_1Bs)\le\ell(W_1U_1)+\ell(Bs)=\ell(W_1U_1)+1\le\ell(W_1)+\ell(U_1)+1=
\ell(W)+\ell(U_1)=\ell(Wt)+\ell(U)-1 = \ell(yVx)+\ell(U)-1 = \ell(U)+\ell(yV)-1 = \ell(Y)-1=n-1$.

\medskip
Case 2.2. $\sup\bold d(Y)<\sup Y$. Recall that $Yx=Y_1Ax=Y_1B$ where $A=\varphi(Y)$ and $B=Ax=\iotaRev(Yx)$.
So, we have $\bold d(Y)=AY_1$. Thus the condition $\sup\bold d(Y)<\sup Y$ reads as
$\sup AY_1 < \sup Y_1A = \sup A + \sup Y_1$, hence, by Lemma \lemSup, we have $AC\preccurlyeq\delta$ where we set
$C=\phiRev(Y_1)$, $Y_1=CY_2$. Since $B=\iotaRev(Yx)$ and $Yx=Y_1B$, we have
$\phiRev(Yx)=\phiRev(Y_1B)=\phiRev(Y_1)=C$.
Thus $(x_1,\bold d(Y))\in\Cal Q_{n-1}$ because
$x_1\bold d(Y)=x_1\bold d(Y_1A)=x_1AY_1=AxY_1\sim Y_1Ax=Yx\sim X$ and
$\bold d(Y)=\bold d(Y_1A)=AY_1=ACY_2$ whence $\ell(\bold d(Y))\le\ell(AC)+\ell(Y_2)=1+\ell(Y_2)=\ell(Y_1)=\ell(Y)-1$.

\medskip
Case 2.3. $\sup_s Y < \sup\bold d(Y)=\sup Y$. Let us show that this case is impossible.
Indeed, the condition $\sup\bold d(Y)=\sup Y$ combined with
Lemma \lemD b yields $\ell_1(Y_1)=0$. Since, moreover, $Yx=Y_1B$, $B=Ax=\iotaRev(Yx)$ and $\|B\|=2$,
we obtain $\ell_1(Yx)=0$. By (\eqCaseII) this implies that $Yx$ is rigid which contradicts (\eqCaseIIa).
\qed\enddemo


\proclaim{ Lemma \lemCaseIII }
Let $(x,Y)\in\Cal Q_{n,p}$, $p<0$.
Suppose that $\ell(x Y)=\ell(Y x) = 1+\ell(Y)$.
Then either $x Y\in\SSS(X)$, or $Y x\in\SSS(X)$, or $\Cal Q_{n-1}\ne\varnothing$, or
 $\Cal Q_{n,p+1}\ne\varnothing$
\endproclaim

\demo{ Proof }
The condition $\ell(Yx)=\ell(Y)+1$ implies
$\varphi(Yx)=x$ and hence $\bold d(Yx)=xY$.

\medskip
Case 1. $\sup Yx > \sup_s Yx$.
By [\refBKLtwo] we then have
$$
      \sup\bold d(xY) = \sup\bold d^2(Yx) < \sup Yx.                      \eqno(\eqCaseIII)
$$
Since $\ell(xY)=\ell(Yx)$, we have $\sup\bold d(Yx)=\sup xY = \sup Yx$. Hence
$\ell_1(Y)=0$ by Lemma \lemD b. Let $A=\varphi(Y)$, $B=\varphi(xY)$, $C=\iota(xY)$, and let
$xY=CUB$.
We have $A\ne B$ (otherwise we would obtain $\sup\bold d(Yx) < \sup Yx$ by Lemma \lemD a)
and we have $A\succcurlyeq B$ by Lemma \lemCharney.
Hence $\|B\|=1$. By combining this fact with $\ell_1(CUB)=\ell_1(xY)=\ell_1(Yx)=1$, we obtain
$\ell_1(CU)=0$. It follows that the left normal form of $\delta^{-p}CU$ coincides with its right normal form,
in particular, $\phiRev(CU)=\iota(CU)=C$.
By (\eqCaseIII), we have $\sup BCU = \sup\bold d(xY) < \sup xY = \sup CUB = \sup B+\sup CU$.
Hence, by Lemma \lemSup, we have $B\phiRev(CU)\preccurlyeq\delta$, that is
$BC\preccurlyeq\delta$. This implies $BC=\delta$ because $\|C\|=2$ (recall that $\ell_1(CU)=0$) and $\|B\|=1$.
We have $x\preccurlyeq\iota(xY)=C$, hence $Bx\preccurlyeq BC=\delta$ which yields $Bx=x_1B$ with
$x_1\in x^G\cap\Cal A$. Since $x\preccurlyeq C$, we may write $C=xC'$, $C'\in[1,\delta]$.
So, for
$Z=BC'U$, we obtain
$x_1 Z= x_1BC'U = BxC'U=BCU = \bold d(CUB)=\bold d(xY)\sim X$ and
$Z=BC'U\sim C'UB=x^{-1} CUB = Y$. We have
$\ell(Z) \le\ell(BC')+\ell(U) = 1+\ell(U) = \ell(xY)-1 = \ell(Y)=n$, hence
$(x_1,Z)\in\Cal Q_{n,p}$. Since $x_1Z=x_1BC'U = BxC'U=BCU=\delta U$, we have
$\inf x_1 Z > \inf U = p$, thus the result follows from Lemma \lemCaseI.

\medskip
Case 2. $\sup Yx = \sup_s Yx$. If $\inf Yx = \inf_s Yx$, then $Yx\in\SSS(X)$ and we are done.
So, we suppose that  $\inf Yx < \inf_s Yx$. Then, by Lemma \lemCDii a, we have
$$
      \inf Yx < \inf\bold c(Yx).                                  \eqno(\eqCaseIIIa)
$$
Let $A=\iota(Y)$, $Y=AY_1$. The condition $\ell(Yx)=\ell(Y)+1$ implies that $\varphi(Y)\cdot x$
is left weighted whence $\iota(Yx)=\iota(Y)=A$.
Thus $\bold c(Y)=Y_1A$, $\bold c(Yx)=Y_1xA$, and
$$
     \varphi(Yx)=\varphi(Y_1x)=x.                                 \eqno(\eqCaseIIIaa)
$$

\medskip
Case 2.1. $\inf Y = \inf_s Y$. Let $t=\frak p(Yx)$, $A=tA'$, thus $\frak s(Yx)=A'Y_1xt$.
Then $xt\preccurlyeq\delta$, hence $xt=tx_2$, $x_2\in x^G\cap\Cal A$, and we obtain
$\frak s(Yx)=Y^tx_2$. By (\eqCaseIIIa) combined with [\refGG; Lemma 4] we have
$$
             \inf Yx < \inf\frak s(Yx).                                   \eqno(\eqCaseIIIb)
$$
Since $t\preccurlyeq A=\iota(Y)$, we have $\ell(Y^t)\le\ell(Y)+1$. If $\ell(Y^t)\le\ell(Y)$, then
the result follows from Lemma \lemCaseI\ applied to $(x_2,Y^t)$, because
$x_2Y^t\sim Y^tx_2=\frak s(Yx)\sim X$
and $\inf Y^t x_2 > p$ by (\eqCaseIIIb). So, we assume that
$$
                 \ell(Y^t) = \ell(Y) + 1.                        \eqno(\eqCaseIIIc)
$$
The condition $t\preccurlyeq A=\iota(Y)$ implies $\inf Y^t \ge\inf Y$. Since $\inf Y=\inf_s Y$, it
follows that $\inf Y^t = \inf_s Y$. Hence, by the `right-to-left' version of
Lemma \lemUxVi, we have $Y^t = UyV$ with $UV\sim z$, $\ell(U)+\ell(V)+1=\ell(Y^t)$, $\ell(V)\ge 1$,
and $\iotaRev(Uy)\cdot\phiRev(V)$ right weighted.
The last two conditions imply $\iotaRev(Y^t)=\iotaRev(V)$; we denote this element by $B$ and we set
$V=V_1 B$.
By (\eqCaseIIIb) and (\eqCaseIIIc) we have
$\inf Y^t + \inf x_2 = \inf Y^t = \inf Y = \inf Yx < \inf\frak s(Yx) = \inf Y^t x_2$, hence
$\delta\preccurlyeq\iotaRev(Y^t)x_2=Bx_2$ by Lemma \lemInf.
Since $\|Bx_2\|\le\|\delta\|$, this means that $Bx_2=\delta$, and we obtain
$\frak s(Yx) = UyVx_2 = UyV_1Bx_2=UyV_1\delta\sim yZ$ where $Z=V_1\delta U$.
Since $UV\sim z$, we have $Z\sim UV_1\delta=UV_1Bx_2=UVx_2\in z^G x^G$. Since, moreover,
$\ell(Z)\le\ell(U)+\ell(V_1)=\ell(U)+\ell(V)-1=\ell(Y^t)-2 = \ell(Y)-1$,
we conclude that
$(y,Z)\in\Cal Q_{n-1}$.

\medskip
Case 2.2. $\inf Y < \inf\bold c(Y)$. Recall that $Y=AY_1$ and $A=\iota(Y)=\iota(Yx)$.
Let  $B=\iotaRev(Y_1)$, $Y_1=Y_2B$.
Since $\inf Y_1+\inf A=\inf Y < \inf\bold c(Y)=\inf Y_1A$, we have
$\delta\preccurlyeq BA$ by Lemma \lemInf.
Hence $B=CD$ and $DA=\delta$ for some simple elements $C$ and $D$.
By Theorem \thQPone, the left normal form of $DxA$ is $D'\cdot x_1\cdot A'$
with $A',D'\in\Cal P$, $x_1\in x^G\cap\Cal A$, and $D'A'=\delta$.

Since $\iota(Yx)=\iota(Y)=A$, we have
$\bold c(Yx)=Y_1xA = (Y_2C)(DxA)$. Hence,
$\delta\preccurlyeq Y_2C\,\iota(DxA)= Y_2CD'$
by (\eqCaseIIIa) combined with Lemma \lemInf.
Hence, for $Z=A'Y_2CD'$, we have $\inf Z > \inf Y$ and $\ell(Z)\le\ell(A')+\ell(Y_2)+\ell(C)+\ell(D')-1
\le\ell(Y)$. Since
$Z\sim Y_2CD'A'=Y_2C\delta=Y_2CDA=Y_1A\sim Y$ and
$x_1 Z\sim Y_2CD'x_1A'=Y_2CDxA=Y_1xA\sim Yx\sim X$, we conclude that $(x_1,Z)\in\Cal Q_{n,p+1}$.

\smallskip
Case 2.3. $\inf Y = \inf\bold c(Y) < \inf_s Y$. Then $\ell_2(Y_1)=0$ by Lemma \lemC b.
By (\eqCaseIIIaa), this implies $\ell_2(Y_1x)=0$ whence $\iotaRev(Y_1x)=\varphi(Y_1x)=x$.
By (\eqCaseIIIa), we have $\inf Y_1xA =\inf\bold c(Yx) > \inf Yx = \inf Y_1x + \inf A$.
Hence $\delta\preccurlyeq\iotaRev(Y_1x)A=xA$ by Lemma \lemInf.
Since $\|xA\|\le 3$, this means that $xA=\delta$. Hence $xA=Ax_1$, $x_1\in\Cal A$,
and we obtain $\bold c(Yx) = Zx_1$ where $Z=Y_1A=\bold c(Y)\sim Y$ and $\delta\preccurlyeq Zx_1$,
so, the result follows from Lemma \lemCaseI\ applied to $(x_1,Z)$.
\qed\enddemo


\proclaim{ Lemma \lemFinal } Let $(x,Y)\in\Cal Q$ and $Yx\in\SSS(X)$.
Then there exist $(x_1,Y_1)\in\Cal Q$ such that $x_1Y_1\in\SSS(X)$.
\endproclaim

\demo{ Proof } Let $A=\iotaRev(Yx)$ and $Yx=UA$. Then $A\succcurlyeq x$ whence $A=sx=x_1s$ and $Y=Us$
for a simple element $s$.
Let $X_1=\cRev(Yx)$ and $Y_1=sU$. Then we have
$X_1=AU=x_1sU=x_1Y_1$, hence $(x_1,Y_1)\in\Cal Q$ and
$x_1Y_1\in\SSS(X)$.
\qed\enddemo

Theorem \thMain\ immediately follows from Lemmas \lemCaseI\ -- \lemFinal.



\head\sectLee. Structure of $\SSS(X)$ when $\|\Delta\|=3$ (after S.-J.~Lee)
\endhead
Here we give a summary of those results from [\refLee; Chapter 4] which extend
to any homogeneous Garside group with Garside element of letter length 3.

Let $(G,\Cal P,\Delta)$ be a homogeneous Garside structure with set of atoms $\Cal A$
such that $\|\Delta\|=3$.

We say that $X\in G$ is {\bf rigid} if $\varphi(X)\cdot\iota(X)$ is left weighted.
Following [\refLee], we say that $X$ is {\bf strictly rigid} if it is rigid and
$\ell_1(X)=0$ or $\ell_2(X)=0$ (see (\eqDefCount)).
If $X\in\USS(X)$, then we define the cycling orbit of $X$ as $O_X=\{\bold c^m\tau^k(X)\mid k,m\ge 0\}$.

\proclaim{ Proposition \propLee } Let $X\in\USS(X)$, $\ell(X)\ge2$. Then:

\smallskip\noindent
(a). $\SC(X)=\USS(X)$. 

\smallskip\noindent
(b). $\SSS(X) = \bigcup_{m\ge 0}{\cRev}^m( \USS(X) )$. 

\smallskip\noindent
(c). One and only one of the following alternatives holds:
\roster
\item "(i)" each element of $\USS(X)$ is strictly rigid and $\SSS(X)=\USS(X)$;
\item "(ii)" each element of $\USS(X)$ is rigid but not strictly rigid, and $\USS(X)=O_X$;
\item "(iii)" no element of $\SSS(X)$ is rigid and $\SSS(X)=\USS(X)=O_X$.
\endroster
\endproclaim

\proclaim{ Lemma \lemLeeCC }
If $X$ is not rigid and $X\in\SSS(X)$, then
$\cRev(\bold c(X))=\dRev(\bold d(X))=X$.
\endproclaim

\demo{ Proof }
 If $\ell(X)=1$, the statement is evident. Assume that $\ell(X)>1$.
Since $X$ is not rigid, the product $\varphi(X)\cdot\iota(X)$ is not left weighted.
Since $X\in\SSS(X)$, this implies
$\|\varphi(X)\|=1$ and $\|\iota(X)\|=2$.
Let $X=\iota(X)U$. Then $\bold c(X)=U\iota(X)$,
hence $\iotaRev(\bold c(X))\succcurlyeq\iota(X)$. Thus fact combined with $\|\iota(X)\|=2$
implies $\iotaRev(\bold c(X))=\iota(X)$ whence $\cRev(\bold c(X))=X$.
Similarly $\dRev(\bold d(X))=X$.
\qed\enddemo

\proclaim{ Lemma \lemNonRigid } Let $X\in G$, $\ell(X)>1$.
Suppose that $X^G$ does not contain any rigid element.
Then $\SC(X)=\SCRev(X)=\SSS(X)$.
\endproclaim

\demo{ Proof } Lemma \lemLeeCC\ implies that $\bold c$ and
$\bold d$ are bijective mappings from $\SSS(X)$ to itself and that $\cRev$ and
$\dRev$ are their inverse mappings. Hence $\frak s$ and $\sRev$ also are bijective mappings
from $\SSS(X)$ to itself.
\qed\enddemo

\demo{ Proof of Proposition \propLee }
(a). 
If $X^G$ does not contain a rigid element, then the result follows from Lemma \lemNonRigid.
Otherwise it follows from [\refBGGMi; Theorem 3.15] which states that if $X^G$ contains a rigid element,
then all elements of $\USS(X)$ are rigid.

\medskip
(b). 
Let $X\in\SSS(X)$ and let $m\ge0$ be the minimal number such that $Y=\bold c^m(X)\in\USS(X)$.
Then 
${\cRev}^m(Y)=X$ by Lemma \lemLeeCC.

\medskip
(c). The fact that $\USS(X)=O_X$ when $X$ is not strictly rigid is proven in [\refLee; Theorem 4.4.1].
All the other statements follow from (a), and [\refBGGMi; Theorem 3.15].
\qed\enddemo


\Refs
\def\r{\ref}

\r\no\refBessis
\by D.~Bessis
\paper The dual braid monoid
\jour Ann. Sci. \'Ecole Norm. Sup. \vol 36 \yr 2003 \issue 5 \pages 647--683
\endref

\r\no\refBessisCorran
\by D.~Bessis, R.~Corran
\paper Non-crossing partitions of type $(e,e,r)$
\jour Adv. Math. \vol 202 \yr 2006 \pages 1--49
\endref

\r\no\refBGGMi
\by  J.~S.~Birman, V.~Gebhardt, J.~Gonz\'alez-Meneses
\paper Conjugacy in Garside groups I: cycling, powers, and rigidity
\jour  Groups, Geom. and Dynamics \vol 1 \yr 2007 \pages 221--279
\endref


\r\no\refBKLone
\by J.~Birman, K.-H.~Ko, S.-J.~Lee
\paper A new approach to the word and conjugacy problems in the braid groups
\jour  Adv. Math. \vol 139 \yr 1998 \pages 322--353
\endref

\r\no\refBKLtwo
\by J.~Birman, K.-H.~Ko, S.-J.~Lee
\paper The infimum, supremum, and geodesic length of a braid conjugacy class
\jour Adv. Math. \vol 164 \yr 2001 \pages 41--56
\endref

\r\no\refBO
\by M.~Boileau, S.~Yu.~Orevkov
Quasipositivit\'e d'une courbe analytique dans une boule pseudo-convexe,
C. R. Acad. Sci. Paris. S\'er. I \vol 332 \yr 2001 \pages 825--830
\endref

\r\no\refCW
\by M.~Calvez, B.~Wiest
\paper A fast solution to the conjugacy problem in the four-strand braid group
\jour J. Group Theory \vol 17 \yr 2014 \pages 757--780
\endref

\r\no\refCharney
\by R.~Charney
\paper Artin groups of finite type are biautomatic
\jour Math. Ann. \vol 292 \yr 1992 \pages 671--683
\endref

\r\no\refDehornoy
\by P.~Dehornoy
\paper Groupes de Garside
\jour Ann. Sci. \'Ecole Norm. Sup. \vol 35 \yr 2002 \pages 267--306
\endref

\r\no\refDP
\by P.~Dehornoy, L.~Paris
\paper Gaussian groups and Garside groups, two generalizations of Artin Groups
\jour  Proc. London Math. Soc. (3) \vol 79 \yr 1999 \pages 569--604
\endref

\r\no\refEM
\by E.~ElRifai, H.~Morton
\paper Algorithms for positive braids
\jour  Quart. J. Math. Oxford Ser. (2) \vol 45 \yr 1994 \pages 479--497
\endref

\r\no \refGarside 
\by     F.~A.~Garside 
\paper  The braid group and other groups 
\jour   Quart. J. Math. \vol 20 \yr 1969 \pages 235--254 
\endref 


\r\no\refGebGM
\by V.~Gebhardt, J.~Gonz\'alez-Meneses
\paper The cyclic sliding operation in Garside groups
\jour  Math. Z. \vol 265 \yr 2010 \pages 85--114
\endref

\r\no\refLee
\by S.-J.~Lee
\paper  Algorithmic solutions to decision problems in the braid groups
\jour Ph.~D. Thesis, 1999
\endref

\r\no\refOrevkovTop
\by S.~Yu.~Orevkov
\paper Link theory and oval arrangements of real algebraic curves
\jour Topology \vol 38 \yr 1999 \pages 779--810
\endref

\r\no\refOrevkovUR
\by S.~Yu.~Orevkov
\paper Quasipositivity test via unitary representations of braid
groups and its applications to real algebraic curves
\jour J. Knot Theory and Ramifications \vol 10 \yr 2001 \pages 1005--1023
\endref

\r\no\refOrevkovGAFA
\by S.~Yu.~Orevkov
\paper Classification of flexible $M$-curves of degree $8$ up to isotopy
\jour GAFA - Geom. and Funct. Anal. \vol 12 \yr 2002 \pages 723--755
\endref

\r\no\refOrevkovQPthree
\by S.~Yu.~Orevkov
\paper Quasipositivity problem for 3-braids
\jour Turkish J. Math. \vol 28 \yr 2004 \pages 89--93
\endref

\r\no\refOrW
\by S.~Yu.~Orevkov
\paper Some examples of real algebraic and real pseudoholomorphic curves
\inbook in: Perspectives in Analysis, Geometry and Topology
\bookinfo Progr. in Math. 296 \publ Birkh\"auser/Springer \publaddr N. Y.
\yr 2012 \pages 355-387 
\endref

\r\no\refQPtwo
\by     S.Yu.~Orevkov
\paper  Algorithmic recognition of quasipositive braids of algebraic length two
\jour   J. of Algebra \vol 423 \yr 2015 \pages 1080--1108 \endref

\r\no\refPicant
\by     M.~Picantin
\paper  Explicit presentations for the dual braid monoids
\jour   C. R. Acad. Sci. Paris. S\'er I \vol 334 \yr 2002 \pages 843-–848
\endref

\r\no\refRudolph
\by L.~Rudolph
\paper Algebraic functions and closed braids
\jour Topology \vol 22 \yr 1983 \pages 191--202
\endref

\endRefs
\enddocument

\enddocument